\documentclass{svjour3}
\usepackage{amssymb,amsfonts}
\usepackage[ruled,norelsize,linesnumbered]{algorithm2e}
\usepackage{graphicx}
\usepackage{textcomp}
\usepackage[dvipsnames]{xcolor}
\usepackage{enumerate}
\usepackage{mathtools}
\usepackage[hidelinks]{hyperref}
\usepackage{silence}
\usepackage{subcaption}

\usepackage[numbers,sort]{natbib}

\usepackage{lmodern}
\usepackage{doi}

\newcommand{\R}{\mathbb{R}}
\newcommand{\N}{\mathbb{N}}

\newcommand{\new}[1]{#1}
\newcommand{\newer}[1]{#1}

\newtheorem{assumption}{Assumption}
\newtheorem{application}{Application}

\DeclareMathOperator*{\argmin}{arg\,min}

\newcommand{\col}{\textnormal{col}}
\newcommand{\T}{\mathcal{T}}
\newcommand{\X}{\mathcal{X}}
\newcommand{\norm}[1]{\bigl|#1\bigr|}

\usepackage{mathtools}
\mathtoolsset{showonlyrefs=true}

\def\BibTeX{{\rm B\kern-.05em{\sc i\kern-.025em b}\kern-.08em
    T\kern-.1667em\lower.7ex\hbox{E}\kern-.125emX}}
\begin{document}

\title{Distributed Nonconvex Optimization with Exponential Convergence Rate via Hybrid Systems Methods}

\author{Katherine R. Hendrickson, Dawn M. Hustig-Schultz, 
Matthew T. Hale, Ricardo G. Sanfelice}
\institute{
Katherine Hendrickson \at
Merlin Labs \\
Boston, MA, USA 02111 \\
kat.hendrickson@merlinlabs.com
\and
Dawn M. Hustig-Schultz \at
University of California, Santa Cruz \\
Santa Cruz, CA, USA 95064 \\
dhustigs@ucsc.edu
\and
Matthew T. Hale, Corresponding Author \at
Georgia Institute of Technology \\
Atlanta, GA, USA 30332 \\
matthale@gatech.edu
\and
Ricardo G. Sanfelice \at
University of California, Santa Cruz \\
Santa Cruz, CA, USA 95064 \\
ricardo@ucsc.edu
}

\authorrunning{K. Hendrickson, D.M. Hustig-Schultz, M.T. Hale, and R.G. Sanfelice}
\titlerunning{Distributed Nonconvex Optimization via Hybrid Systems Methods}


\maketitle

\begin{abstract}
We present a hybrid systems framework for distributed multi-agent optimization in which
agents execute computations in continuous time and communicate in discrete
time. The optimization algorithm is 
analogous to 
a continuous-time form of parallelized
coordinate descent. Agents implement an update-and-hold strategy
in which gradients are computed at communication times and held
constant during flows between communications.
The completeness of solutions under these hybrid dynamics is established. 
Then, we prove that this system is globally
exponentially stable to a minimizer of a possibly nonconvex, smooth objective function that satisfies the Polyak-{\L}ojasiewicz (PL) condition. Simulation results are presented for three different applications and illustrate the convergence rates
and the impact of initial conditions upon convergence. 
\end{abstract} \\



\section{Introduction}
\subsection{Motivation}
Optimization problems arise in many
areas of engineering, including machine learning~\cite{sra12},
communications~\cite{luo06}, robotics~\cite{verscheure09}, and others. 
Across all application areas, the goal is to design
an algorithm that will converge to a minimum of an objective function,
possibly under some constraints. Recently, 
there has been increased interest in studying 
optimization algorithms
in continuous time to use
tools from dynamical systems to establish convergence
to minimizers, e.g., in~\cite{su16,rahili2017,garg20}. While a large body of optimization work focuses on convex optimization, nonconvex problems often arise in a variety of fields, including machine learning~\cite{du2019gradient, fazel2018global, lian2015asynchronous} and communication networks~\cite{Chiang2009}, and there has arisen interest in establishing convergence guarantees for non-convex problems. 

In this paper, 
we develop a multi-agent framework for nonconvex optimization in which agents' computations are modeled in continuous time. This is motivated by two factors. First, we wish to leverage the large collection of 
tools from dynamical systems to analyze multi-agent optimization and connect to the
growing body of work that uses continuous-time models of computation. 
Second, there also exist controllers for multi-agent systems
that are designed to operate in continuous time to minimize some objective function,
e.g., in consensus~\cite{olfati07} and coverage control~\cite{cortes04},
and our analyses will connect our work to such systems. 
However, while
individual agents' computations occur in continuous time, their
communications are most naturally modeled in discrete time 
because communicated information arrives at its recipients at
individual instants in time. Thus, the joint modeling and analysis
of agents' computations and communications creates a mixture
of continuous- and discrete-time elements, which leads us
to develop a hybrid systems framework for multi-agent optimization. 

The framework that we develop can solve a class of problems that includes some non-convex problems.
In particular, we consider smooth objective functions that satisfy the Polyak-{\L}ojasiewicz (PL) inequality~\cite{polyak1963gradient}. Recent interest in the PL inequality and related properties has led to the development of discrete-time nonconvex optimization approaches~\cite{karimi2016linear, charles18,loizoustoch20, Gower2021SGDFS, Bassily2018OnEC}, including distributed algorithms~\cite{yazdani2021asynchronous,ubl22,ubl23}. Problems that satisfy the PL inequality include matrix factorization~\cite{sunmatrix15}, minimizing logistic loss over a compact set~\cite{karimi2016linear}, and the training of some neural networks~\cite{charles18}.
The set of functions that satisfy the PL inequality also includes those that are strongly convex,
and our developments therefore apply to strongly convex functions, which are commonly studied in distributed
optimization settings~\cite{behrendt23}. 

The algorithm that we develop is analogous to a hybrid systems version 
of parallelized block coordinate descent~\cite{bertsekas89}, in which each agent
updates only a subset of all decision variables 
using continuous-time computations
and agents
communicate these updated values in discrete time to other agents. 
In the framework that we develop, all agents' communications are intermittent; 
agents' communications
occur when a timer reaches zero, at which point the timer
is reset to some value within a specified range. 
Agents use a sample-and-hold strategy in which gradients
are computed at the communication times and then held constant
and used in computations until the next communication event. 
This approach is inspired by recent work~\cite{PHILLIPS2019} that has
successfully applied it to synchronization problems.

\subsection{Contributions}
We leverage analytical tools from the theory of hybrid systems to prove that
this algorithmic framework has several desirable properties, and
our contributions are: 
\begin{itemize}
    \item We define a hybrid system model for distributed optimization. To the best of our knowledge, this is the first distributed hybrid system model that uses a parallelized approach for nonconvex problems rather than a consensus-based approach.
    \item We show that under our model, every maximal solution is complete, with the time domain allowing arbitrarily large ordinary time~$t$. As a result, there are no theoretical obstructions to running this algorithm for arbitrarily long periods of time.
    \item We use Lyapunov analysis to show that, even under intermittent information sharing, the hybrid optimization algorithm is globally exponentially stable to a minimizer of an objective function, and we derive an explicit convergence rate in terms of system parameters. 
    \item We show robustness to inaccuracies in the measurement of the times at which communication events occur. 
    \item Finally, we present three different applications, including those with nonconvex objective functions, that demonstrate the performance of our model.
\end{itemize}

\subsection{Related Work}
The developments in this paper can be regarded as hybrid counterparts to
``classical'' discrete-time algorithms in multi-agent 
optimization~\cite{bertsekas89}. 
Related research in multi-agent continuous-time optimization includes~\cite{gharesifard2014,lu2012,rahili2017}, though 
those works all use a consensus-based optimization framework in which
both computations and communications are modeled as occurring in continuous time. 
In this work, we avoid continuous-time
communications in order to model problems in which constant communications are not possible, e.g.,
over long distances, or simply undesirable, e.g., when battery power is limited. 

The closest works to the current paper are~\cite{KIA2015},~\cite{PHILLIPS2019}, 
which also study multi-agent optimization algorithms with continuous-time computations and
discrete-time communications. However, those works also use consensus-based optimization
algorithms in which each agent has a local objective function, computes
new values for all decision variables, and averages its decision variables with other agents. 
In contrast, we consider all agents having a common objective function and 
we only require each agent to compute updates to a small subset of the decision
variables in a problem. 
\new{
This approach has the benefit that each agent's computational burden grows slowly as 
a problem grows
since each agent updates only~$\frac{d}{N}$ decision variables on average. 
When~$N$ is large, the value of~$\frac{d}{N}$ can grow quite slowly as a function of~$d$, which  
results in only small increases in each agent's computational burden. 
}

In addition, the hybrid system model that we develop offers several analytical features.
First, existing block coordinate descent algorithms are typically modeled in discrete time.
When used in a continuous-time system, these types of
discrete-time computations will be done with samples of continuous-time 
state values, though convergence analyses for the discrete-time updates 
will apply only to the samples, not to the 
continuously evolving intersample state values. 
However, within the
hybrid framework in this paper, we analyze the time evolution of both the sampled state values
and the intersample state values, which characterizes state evolution
at all points in time. 
Second, when a hybrid system is well-posed and 
has a compact pre-asymptotically stable set,  
it follows that such pre-asymptotic stability is robust to small perturbations~\cite[Theorem~7.21]{goebel12}. 
In this paper, we show that this robustness applies when solving the problems we consider, and
thus our use of a hybrid model lets our analysis automatically inherit these robustness properties. 

This paper is an extension of the conference paper~\cite{hendricksonhybrid21} 
which applied only to strongly convex objective functions. This paper modifies 
the previous hybrid system model to provide global convergence and reformulates 
all previous results to apply to objective functions that satisfy the 
PL inequality. 
\new{
This class of functions includes some non-convex functions, and all 
theoretical results and proofs from~\cite{hendricksonhybrid21} have been 
reformulated to accommodate this nonconvexity.
In particular, generic PL functions can have any number of local minima
and need not have strongly monotone gradients, both of which differ
from strongly convex functions. 
As described above, 
we are motivated to analyze PL functions because 
they appear in a wide range of engineering problems,
and this paper therefore extends our hybrid algorithm to apply to those problems. 
}
Moreover, a tighter convergence rate is established and more 
applications are demonstrated via simulation. Finally, new results regarding 
robustness to perturbations are presented.

\subsection{Organization}
The rest of the paper is organized as follows. Section~\ref{sec:prob} includes our problem statement, assumptions, and algorithm. We then present our hybrid system model for multi-agent optimization 
in Section~\ref{sec:hsm} and establish the existence of complete solutions. Section~\ref{sec:stability} proves
that the hybrid multi-agent update law is globally 
exponentially stable, and then Section~\ref{sec:robustness} shows that this exponential
stability guarantee is robust to a certain class of perturbations. 
We include numerical results in Section~\ref{sec:sim}, and we present our conclusions in Section~\ref{sec:concl}. 

\section{Problem Statement and Algorithm Overview} \label{sec:prob}
In this section, we state the class of problems that we consider, and
we give an overview of the hybrid optimization algorithm that is the focus
of the remainder of the paper. First, we present some general notation.

\textbf{Notation and Terminology.} 
Let~$\R$ denote the set of real numbers,
let~$\R_{\geq 0}$ denote the non-negative reals,
and let~$\R_{> 0}$ denote the positive reals. 
Let~$\mathbb{N}$ denote the non-negative integers, and
let~$\mathbb{N}_{>0}$ denote the positive integers. 
For~$d \in \mathbb{N}_{> 0}$, 
let~$\boldsymbol{0}_d$ be a vector of zeros in~$\R^d$ 
and~$\boldsymbol{1}_d$ be a vector of ones in~$\R^d$. Define the set~$[p]:=\lbrace 1, 2, \ldots , p\rbrace$ for 
any~$p \in \N_{>0}$. 
For vectors~$x_1, x_2, \dots, x_n$, define~$\col (x_1, x_2, \dots, x_n) := (x_1^{\top}, x_2^{\top}, \dots, x_n^{\top})^{\top}$. Throughout the paper,~$\norm{\cdot}$ denotes the Euclidean norm.
We use~$\textnormal{dom} f$ to denote the domain of a function~$f$. 
The set~$\mathcal{K}_{\infty}$ denotes the set of 
class~$\mathcal{K}_{\infty}$ functions, i.e.,
functions~$\alpha : \R_{\geq 0} \to \R_{\geq 0}$ that are (i) strictly increasing,
(ii) satisfy~$\alpha(0) = 0$, and (iii) satisfy
$\lim_{r \to \infty} \alpha(r) = \infty$. 
The set~$\mathcal{KL}$ denotes the set
of class-$\mathcal{KL}$ functions, i.e., 
functions~$\gamma : \mathbb{R}_{\geq 0} \times \mathbb{R}_{\geq 0} \to \mathbb{R}_{\geq 0}$ 
such that (i) $\gamma$ is non-decreasing in its first argument, (ii) $\gamma$
is non-increasing in its
second argument, (iii)
$\lim_{r \to 0^+} \gamma(r, s) = 0$
for each~$s \in \R_{\geq 0}$, and
(iv) $\lim_{s \to \infty} \gamma(r, s) = 0$ for each~$r \in \R_{\geq 0}$. 
We use ``ODE'' to mean ``ordinary differential equation''. 

\subsection{Problem Formulation}
We consider a group of agents jointly solving an optimization problem that may be nonconvex. Suppose there are~$N$ agents
that will each execute computations locally and then share the results of those computations with other agents. 
For scalability, only a single agent will update each decision variable. 
In many practical settings, we expect bandwidth to be limited and/or agents to
have limited onboard power available, which means communications should not be constant.

Under these conditions, we consider minimization problems of the following form: 
\begin{problem} \label{prob:first}
Given an objective function~$L : \R^n \to \R$,
\begin{equation}
\textnormal{minimize } L(x), \qquad x \in \R^n
\end{equation}
while requiring that (i) only one agent updates any entry of the decision variable~$x$, and (ii) agents
require only sporadic information sharing from others. 
\end{problem}

We first assume the following about the objective function~$L$. 

\begin{assumption} \label{as:L}
The function~$L$ is twice continuously
differentiable and~$K$-smooth (namely, $\nabla L$ is~$K$-Lipschitz). 
\hfill $\triangle$
\end{assumption}

Rather than requiring that~$L$ be convex, we will instead consider a class functions that includes some nonconvex problems. In particular, we are interested in problems that satisfy the Polyak-{\L}ojasiewicz (PL) inequality~\cite{karimi2016linear}.

\begin{assumption}\label{as:pl}
The set of stationary points of~$L$, defined as
$\X^* = \{x \in \mathbb{R}^n : \nabla L(x) = 0\}$,
is non-empty, and
the function $L$ satisfies the Polyak-{\L}ojasiewicz (PL) inequality. Namely, there exists some constant~$\beta > 0$ such that
\begin{equation}
    \frac{1}{2}\norm{\nabla L(x)}^{2} \geq \beta \bigl(L(x)- L(x^*) \bigr)  
\end{equation}
for all~$x \in \R^n$ and 
all~$x^* \in \X^*$. \hfill $\triangle$ 
\end{assumption}

If a function satisfies the PL inequality with constant~$\beta$, we say that the function is~$\beta$-PL. One useful property resulting from Assumption~\ref{as:pl} is that all stationary points of~$L$ are global minima. 
\begin{lemma}
Let Assumption~\ref{as:pl} hold. Then every local minimizer of~$L$ is a global minimizer.
\end{lemma}
\noindent \emph{Proof: }
See Section~2.2 of~\cite{karimi2016linear}. 
\hfill $\qed$

We define
\begin{equation} \label{eq:lstar}
L^* := L(x^*)
\end{equation}
as the global minimum value of~$L$, and,
under Assumption~\ref{as:pl} it is attained at 
every~$x^* \in \X^*$. 
Assumptions~\ref{as:L} and~\ref{as:pl} ensure that even nonconvex problems still retain some geometric structure that allows for convergence analysis. In particular, the combination of these two assumptions provides both an upper and lower bound on a function's gradient, which will be useful in the forthcoming analysis. 
\begin{lemma}
For a function~$L$ that satisfies Assumptions~\ref{as:L} and~\ref{as:pl},
for all~$x \in \R^n$ and~$x^* \in \X^*$ we have
    $2 \beta \left(L(x) - L^* \right) \leq \norm{\nabla L(x)}^{2} 
    \leq K^2 \norm{x - x^*}^2$, 
where~$L^*$ is from~\eqref{eq:lstar},
$K$ is the Lipschitz constant of~$\nabla L$ from Assumption~\ref{as:L},
and~$\beta$ is the PL-constant of~$L$ from Assumption~\ref{as:pl}. 
\end{lemma}
\emph{Proof:}
The left inequality follows directly
from Assumption~\ref{as:pl}. The right inequality follows
by noting that~$\nabla L(x^*) = 0$
and that therefore~$\norm{\nabla L(x)} = \norm{\nabla L(x) - \nabla L(x^*)}$, 
and then applying the Lipschitz property
of~$L$ from Assumption~\ref{as:L}. 
\hfill $\qed$

We refer the reader to~\cite{karimi2016linear} for a thorough discussion of the PL condition in relation to other function properties. Among the strong convexity, essential strong convexity, weak strong convexity, restricted secant inequality, error bound, PL, and quadratic growth conditions, the authors of~\cite{karimi2016linear} establish that the PL and error bound conditions are the most general under which linear convergence to minimizers is achieved. In fact, given our Assumption~\ref{as:L},~\cite[Theorem 2]{karimi2016linear} establishes that any function satisfying the strong convexity, essential strong convexity, weak strong convexity, restricted secant inequality, or error bound conditions also satisfies the PL condition. Thus, we enforce the PL condition as an assumption because it unifies a wide range of problems.

\subsection{Mathematical Framework}
We solve Problem~\ref{prob:first} by developing a hybrid systems framework
in which agents optimize with decentralized gradient descent in continuous time
and communicate their iterates with each other in discrete time. 
Analogously to past research that has developed distributed versions of the discrete-time gradient
descent law, our update law begins with the 
(centralized) 
first-order dynamical system
\begin{equation} \label{eq:gdall}
    \dot{x} + \nabla L(x) =0.
\end{equation}
This is motivated by the use of gradient-based controllers in multi-agent systems, e.g., in consensus~\cite{olfati07}, as well
as the simplicity of distributing gradient-based updates and the robustness 
to intermittency of communications that results from doing so~\cite{bertsekas89}. 

We seek to distribute~$\eqref{eq:gdall}$ across a team of agents in
accordance with the parallelization requirement 
in Problem~\ref{prob:first}. 
We consider~$N$ agents indexed over~$i \in [N]$ and divide~$x \in \R^n$ into~$N$ blocks. Then agent~$i$ is responsible for updating and 
communicating values of the~$i$-th block,~$x_{i} \in \mathbb{R}^{n_i}$, 
where~$n_i \in \mathbb{N}_{> 0}$ and~$\sum_{i \in [N]} n_i = n$. 
Thus, the variable~$x$ may be written as the vertical concatenation of all agents' blocks, i.e.,~$x = \col (x_1, x_2, \dots, x_N)$. Each agent performs gradient descent on their own block 
but does not perform computations on any others. 

Agents' updates occur in continuous time while communications of these updates occur in discrete time. These communication events are coordinated for all agents using a shared timer~$\tau$. When this timer reaches zero, agents will broadcast their current state~$x_i$ to all other agents. The timer will then be reset to a value within a specified interval~$[\tau_{\min}, \tau_{\max}]$. 
Without loss of generality, 
we assume that communications are received at the same time as they are sent (allowing for communication delays requires only adding the length of
delay onto the time between communications). 
When~$\tau = 0$, state values are communicated by agent~$i$ for all~$i \in [N]$ and received by agent~$\ell$ for all~$\ell \in [N]$, and then these communicated values are gathered by agent~$\ell$ into a vector~$\eta^\ell \in \mathbb{R}^n$,
with the received value of~$x_i$ being assigned to~$\eta^\ell_i$. Note that for two agents~$i$ and~$\ell$, the entries~$\eta^i_k$ and~$\eta^\ell_k$ for some~$k \in [N]$ may not be equal at initialization. They will be equal, however, after at least one communication event and will remain equal for the rest of the run of the algorithm.
\new{
We consider the possibility of non-equal initial conditions for two reasons. First, it may be difficult to enforce equality of agents' initial
conditions among many autonomous decision-makers that lack a central coordinator. For example, agents in a large decentralized
network may choose their own initial conditions without coordinating with neighboring agents.
Second, analysis of non-equal initial conditions allows us to establish global convergence results,
where ``global'' indicates that convergence occurs regardless of the initial conditions. 
Globality allows our convergence results to apply to many more scenarios, such as those in which 
agents inadvertently have small disagreements about initial conditions. It also enables us to show
that the hybrid system model we develop is robust to errors in the timing of agents' communications
and computations by employing robust asymptotic
stability tools for hybrid systems (see Section~\ref{sec:robustness}). 
}

The value of~$\eta^i$ is used in agent~$i$'s continuous-time computations in an update-and-hold
manner between communications. Formally, agent~$i$ executes
\begin{equation}
\dot{x}_i = -\nabla_i L(\eta^i), \label{eq:dynamics}
\end{equation}
where the gradient of the function~$L$ with respect to the~$i^{th}$ block and evaluated at some vector~$x$ is written as $\nabla_i L(x) = \frac{\partial}{\partial x_i}L(x)$. This sample-and-hold method is common in the literature~\cite{PHILLIPS2019, phillips2018} and
is used to demonstrate the feasibility of the hybrid approach in multi-agent
optimization. The complete algorithm is summarized in Algorithm~\ref{alg:first}. 

\begin{algorithm}
\SetAlgoLined
Initialization: set~$\eta^i \in \R^n$,~$x_i = \eta^i_i \in \R^{n_i}$, and~$\tau \in [0, \tau_{\max}]$, for all~$i\in[N] $\;
\For{$i \in\lbrace 1, \dots, N \rbrace $}{
\While{$\tau \geq 0$}{
$\dot{x}_i = - \nabla_i L(\eta^i)$\;
$\dot{\tau} = -1$\;
\If{$\tau=0$}{
communicate~$x_i$ to all other agents: reset~$\eta^\ell_i$ to~$x_i$ for all~$\ell \in [N]$\;
reset~$\tau$ to a value in~$ [\tau_{\min},\tau_{\max}]$\;}}
}
\caption{Distributed Gradient Descent}
\label{alg:first}
\end{algorithm}

\new{
The ODE in line~$4$ does not need to be solved in closed form, in the sense of finding
a function of time~$t \to x_i(t)$ that obeys the ODE at all times. 
To implement line~$4$, 
for all~$i \in [N]$, 
the~$i^{th}$ agent only needs to allow~$x_i$ to flow along the direction~$-\nabla_i L(\eta^i)$ 
until the~$\tau = 0$ condition is satisfied in line~$6$. 
Line~$7$ requires agent~$i$ to ``communicate~$x_i$ to all other agents'' and this step
only requires agent~$i$ to communicate its value of~$x_i$ when
the~$\tau = 0$ condition is reached. Agent~$i$ does not need to communicate the time history
of~$x_i$, nor does it even need to store it. 
} The next section develops the hybrid system model that will be
used to analyze Algorithm~\ref{alg:first}. 

\section{Hybrid System Model} \label{sec:hsm}
In this section, we define a hybrid system model that encompasses all agents' current states and their most recently communicated state values. 
Towards defining this ``combined hybrid system'', we first formally state
what constitutes a hybrid system, then we define
the timer that governs communications. This timer allows us to define the hybrid subsystems that are distributed across the agents. Building on that definition, we then present a definition of the combined hybrid system that will be the focus of our analysis, and we verify that this model meets the ``hybrid basic conditions,'' which are defined below. Finally, we show the existence of solutions and conclude that all maximal solutions are complete.

\subsection{Hybrid System Definitions}
For the purposes of this paper, a hybrid system~$\mathcal{H}$ has data~$(C,f,D,G)$ that takes the general form
\begin{align}
    \mathcal{H} = \begin{cases}
    \dot{x} = f(x) & x\in C \\
    x^+ \in G(x) &x\in D
    \end{cases}, \label{eq:hybriddef}
\end{align}
where~$x\in \R^n$ is the system's state,~$f$ defines the flow map and continuous dynamics for which~$C$ is the flow set, 
and~$G$ is the set-valued 
jump map which captures the system's discrete behavior for the jump set~$D$. 
\new{
Here (and below) we use the standard notational convention in which~$x^+$ denotes the value of the
state~$x$ after it undergoes a jump. 
The meaning of~\eqref{eq:hybriddef} is that the state~$x$ evolves according to the 
ODE~$\dot{x} = f(x)$, which defines the
flow dynamics
whenever~$x \in C$, and the state~$x$ undergoes a jump whenever~$x \in D$, which occurs instantaneously,
leading to the difference inclusion~$x^+ \in G(x)$. 
Note that while~$f$ is a single-valued map, $G$ is a set-valued map, meaning that at jumps,~$x$ can take
any value in~$G(x)$. 
}

More information on this definition and hybrid systems can be found in~\cite{goebel12}.

\begin{definition}[Hybrid Basic Conditions,~\cite{goebel12}] \label{def:hybridbasic} A hybrid system~$\mathcal{H} = (C,f,D,G)$ with data given by~\eqref{eq:hybriddef} satisfies the~\emph{hybrid basic conditions} if
\begin{itemize}
    \item~$C$ and~$D$ are closed subsets of~$\R^n$;
    \item~$f$ is defined on~$C$ and is a continuous function from~$C$ to~$\R^n$;
    \item~$G:\R^n \rightrightarrows \R^n$ is outer semicontinuous and locally bounded relative to~$D$, and~$D \subset \textnormal{dom } G$.
\end{itemize}
\end{definition}
If a hybrid system meets the hybrid basic conditions, then we say that the system is~\emph{well-posed} (Theorem 6.30,~\cite{goebel12}).
Well-posedness is desirable because it lets us establish the robustness
of a hybrid system to perturbations, which we do in Section~\ref{sec:robustness}. 

\new{
The following elementary example illustrates the formulation of a hybrid model.}

\begin{example}[From~\cite{ferrante16}]
\new{
Consider the linear time-invariant system
\begin{align}
\dot{z} &= Az + Bu \label{eq:exzfirst}\\
      y &= Mz,
\end{align}
where~$z \in \R^n$ is the state, $y \in \R^q$ is the output,
and~$u \in \R^p$ is the input. 
The input~$u : [0, \infty) \to \R^p$ is measurable
and locally bounded, and~$A \in \R^{n \times n}$,
$B \in \R^{m \times n}$, and~$M \in \R^{q \times n}$
are constant matrices. 
Consider the problem of designing an estimator
for~$z$ that operates when only sporadic measurements
of the output~$y$ are available; see~\cite{ferrante16}. 
Mathematically,
the estimator only has access to outputs~$y(t_k)$
for~$k \in \mathbb{N}_{>0}$ for some collection
of points in time~$\{t_k\}_{k \in \mathbb{N}_{>0}}$.
}

\new{
We assume that the sequence~$\{t_k\}_{k \in \mathbb{N}_{> 0}}$
is strictly increasing and unbounded. We also assume
that there are two constants~$0 < T_1 \leq T_2$ such that
\begin{align}
0 &\leq t_1 \leq T_2 \label{eq:ext1bound} \\
T_1 &\leq t_{k+1} - t_k \leq T_2 \textnormal{ for all } k \in \mathbb{N}_{>0}. 
\end{align}
The value of~$T_1$ is the minimum amount of time that
elapses between~$t_k$ and~$t_{k+1}$ for any~$k \in \N_{> 0}$,
and the value of~$T_2$ is the corresponding maximum amount of time. 
}

\new{
The state estimator generates an estimate~$\hat{z} \in \R^n$.
Given a solution~$t \to z(t)$ to~\eqref{eq:exzfirst} obtained with input~$t \to u(t)$ and resulting
in the output~$t \to y(t)$, 
measurements~$y(t_k)$ occur at discrete time 
instances~$\{t_k\}_{k \in \N_{>0}}$, 
and
the state estimate~$\hat{z}$ evolves in continuous time
between measurements  in order to mirror the continuous-time
dynamics of~\eqref{eq:exzfirst}. 
An estimate of~$t \to z(t)$ is given by~$t \to \hat{z}(t)$ satisfying
\begin{alignat}{2}
&\dot{\hat{z}}(t) = A\hat{z}(t) + Bu(t) \qquad &&\textnormal{ for all } t \neq t_k, k \in \N_{> 0} \\
&\hat{z}(t^+) = \hat{z}(t) + L(y(t) - M\hat{z}(t)) \qquad &&\textnormal{ for all } t = t_k, k \in \N_{> 0},
\end{alignat}
where~$L \in \R^{n \times n}$
and~$\hat{z}(t^+)$ is the right limit of~$t \to \hat{z}(t)$ at~$t = t_k$. 
When~$t \neq t_k$, the state estimate flows in continuous time
according to the system dynamics in~\eqref{eq:exzfirst}. 
When~$t = t_k$, the state estimate jumps instantaneously
from the value~$z(t_k)$ 
to the value~$z(t_k^+)$, which is a correction
term that accounts for any differences between
the measured output~$y(t_k)$ and the predited
output~$M\hat{z}(t_k)$. 
As is common in state estimation problems, 
we examine the estimation
error~$t \to \omega(t) = z(t) - \hat{z}(t)$, where
\begin{alignat}{2}
&\dot{\omega}(t) = A\omega(t) \qquad &&\textnormal{ for all } t \neq t_k, k \in \N_{> 0} \\
&\omega(t^+)  = (I - LM)\omega(t) \qquad &&\textnormal{ for all } t = t_k, k \in \N_{> 0}.  
\end{alignat}
}

\new{
We can formulate the~$\omega$ dynamics as a hybrid system.
The measurement of outputs is driven 
by the sequence~$\{t_k\}_{k \in \N_{> 0}}$, 
and we reformulate
it to be state-driven to 
capture all possible sequences satisfying~\eqref{eq:ext1bound}. 
We introduce a new state~$\tau$ for this purpose, 
and~$\tau$ evolves as follows. Between jumps,
the value of~$\tau$ counts down with unit rate.
When~$\tau = 0$, a jump is triggered and~$\tau$
is reset to a value in~$[T_1, T_2]$ before it
starts counting down again. 
The joint dynamics of~$\omega$ and~$\tau$ can be written
as the hybrid system~$\mathcal{H}_{\omega}$ with dynamics
\begin{equation}
\mathcal{H}_{\omega} 
\left\{
\begin{aligned}
&\left.
\begin{aligned}
\dot{\omega} &= A\omega \\
\dot{\tau} &= -1 
\end{aligned}
\hspace{2cm}
\right\} (\omega, \tau) \in C\\
&\left.
\begin{aligned}
\omega^+ &= (I - LM)\omega \\
\tau^+ &\in [T_1, T_2] 
\end{aligned}
\hspace{0.6cm} 
\right\} (\omega, \tau) \in D
\end{aligned}
\right.,
\end{equation}
where the flow set~$C$ and jump set~$D$ are defined as
\begin{align}
C &= \{(\omega, \tau) \in \R^n \times \R_{\geq 0} : \tau \in [0, T_2]\} \\
D &= \{(\omega, \tau) \in \R^n \times \R_{\geq 0} : \tau  = 0\}.
\end{align}
The sets~$C$ and~$D$ are not disjoint in this model
and, in particular, both~$C$ and~$D$ contain points
of the form~$(\omega, 0)$. 
This property conveys the fact that the timer~$\tau$
counts down from a positive number until it reaches
zero exactly and only then is a jump triggered. 
}

\new{
A complete hybrid model requires the definition
of a flow map~$f$ over~$C$ and a jump map~$G$ 
over~$D$. We define a new state vector~$x := (\omega^T, \tau)^T$. 
Then the flow and jump maps are
\begin{equation}
f(x) := \left(\begin{array}{c}
A\omega \\
-1
\end{array}\right) \qquad \textnormal{ for all } x \in C
\end{equation}
and
\begin{equation}
G(x) := 
\left(\begin{array}{c}
(I - LM)\omega \\
{[T_1, T_2]}
\end{array}\right)  \qquad \textnormal{ for all } x \in D,
\end{equation}
respectively. 
}

\new{
The dynamics of the estimator are entirely encapsulated
in the hybrid model~$\mathcal{H}_{\omega} = (C, f, D, G)$, and the formulation
of this hybrid model enables the use of a large collection
of theoretical tools for its analysis~\cite{goebel12,sanfelice21}. 
}
\hfill $\Diamond$
\end{example}

\new{
Hybrid systems are often analyzed for similar properties as non-hybrid systems, such as asymptotic stability. 
They can also exhibit several types of undesirable behavior. 
For example, hybrid systems can exhibit the Zeno phenomenon in which they undergo an infinite
number of jumps in finite time. They can also exhibit finite escape time in which
some states go to infinity in finite ordinary (flow) time. Both Zeno behavior and finite escape time
cause solutions to have bounded domains. 
In Sections~\ref{ss:comms} through~\ref{ss:combined} we formally define the hybrid system model of 
Algorithm~\ref{alg:first}, and we show in Section~\ref{ss:basicconditions} that it exhibits
neither Zeno behavior nor finite escape time. Then
in Section~\ref{sec:stability} we prove that the states of the hybrid model
asymptotically converge to a minimizer, and we derive a convergence rate. 
}

\subsection{Mechanism Governing the Communication Events} \label{ss:comms}
We seek to account for intermittent communication events that occur only at some time instances~$t_j$, for~$j\in \N_{> 0}$, that are not known \emph{a priori}. We assume that the sequence~$\lbrace t_j\rbrace_{j=1}^\infty$ is strictly increasing and unbounded. Between consecutive time events, some amount of time elapses which we upper and lower bound with positive scalars~$\tau_{\min}$ and~$\tau_{\max}$:
\begin{align}
    0<\tau_{\min} \leq t_{j+1} - t_j \leq \tau_{\max} \quad \textnormal{ for all } j \in \N_{> 0}. \label{eq:timer}
\end{align}
The upper bound~$\tau_{\max}$ prevents infinitely long communication delays and ensures convergence, while the lower bound~$\tau_{\min}$ rules out Zeno behavior.

To generate communication events at times~$t_j$ satisfying~\eqref{eq:timer}, let~$\tau$ be the timer that governs when agents exchange data, where~$\tau$ is defined by
\begin{align}
    \dot{\tau} &= -1   &&\tau \in [0, \tau_{\max}], \label{eq:tau1} \\
    \tau^+ &\in [ \tau_{\min}, \tau_{\max} ] &&\tau = 0, \label{eq:tau2}
\end{align}
for~$\tau_{\min}, \tau_{\max} \in \mathbb{R}_{> 0}$,
\new{where~\eqref{eq:tau1} specifies the flow behavior of~$\tau$ and~\eqref{eq:tau2} specifies
the jump behavior of~$\tau$. 
In words, the timer~$\tau$ flows with~$\dot{\tau} = -1$
until it reaches~$\tau = 0$.
At that point, it stops flowing and a jump is triggered, which resets~$\tau$ 
to a value within~$[ \tau_{\min}, \tau_{\max} ]$, and this new value is denoted by~$\tau^+$.
Since~$\tau^+ \geq \tau_{\min} > 0$, the timer~$\tau$ resumes flowing and this process repeats indefinitely.
}
There is indeterminacy built into the timer in that the reset
map is only confined to a compact interval,~$[\tau_{\min}, \tau_{\max}]$, where~$\tau_{\min}$ and~$\tau_{\max}$ are both positive real numbers.

\subsection{Hybrid Subsystems} 
Recall that for all~$i \in [N]$ agent~$i$ stores their own state variable~$x_i \in \R^{n_i}$. Communications received from all other agents are stored in~$\eta^i \in \R^{n}$, including agent~$i$'s state at the most recent communication event. We define the state of agent~$i$'s hybrid system as~$\xi^i = (x_{i}, \eta^i, \tau)$, where~$x_{i}$ is the state of agent~$i$'s block of the decision variable~$x$ (the one it is responsible for updating),~$\eta^i$ is the vector of state values communicated to agent~$i$ at communication events, and~$\tau$ is defined as in~\eqref{eq:tau1} 
and~\eqref{eq:tau2}. 
Applying the dynamics given in~\eqref{eq:dynamics}, this 
setup leads to the hybrid subsystem 
\begin{align}
    \dot{\xi^i} & = \begin{bmatrix}
    -\nabla_{i} L(\eta^i) \\
    \boldsymbol{0}_n \\
    -1
    \end{bmatrix} &&\xi^i \in \R^{n_i} \times \R^n \times [0,\tau_{\max}] \\
    \xi^{i+} &\in \begin{bmatrix}
    x_{i} \\
    x \\
    [\tau_{\min}, \tau_{\max}]
    \end{bmatrix} &&\xi^i \in \R^{n_i} \times \R^n \times \{0\}.
\end{align}

\new{
The flow and jump sets for~$\xi^i$ have non-empty intersection because
they incorporate the~$\tau$ dynamics from~\eqref{eq:tau1} and~\eqref{eq:tau2} above.
This property ensures that the timer~$\tau$ decreases 
all the way to~$0$ before it resets. 
}
\new{
Subsystem~$i$ has~$\dot{\eta}^i = 0$ for all~$i \in [N]$ to model agents' sample-and-hold strategy.
Agent~$i$ updates only the value of~$x_i$ using the memorized values in~$\eta^i$, 
which is why it has~$\dot{x}_i = -\nabla_i L(\eta^i)$ during flows.
Agent~$i$ does not update~$x_j$ for any~$j \neq i$, and the value of~$x_j$ that agent~$i$ has access to
only changes when agent~$j$ communicates it to agent~$i$. The value of~$x_j$ onboard agent~$i$
is denoted~$\eta^i_j$, and agent~$i$ has the dynamics~$\dot{\eta}^i_j = 0$ to model the fact
that agent~$i$ does not change this value. 
Agent~$i$ does receive the value of~$x_j$ during jumps, which is modeled
by having~$\eta^{i+} = x$ in the jump map. 
Agent~$i$ has~$\dot{\eta}^i_i = 0$ because its computations are modeled
by the dynamics of the state~$x_i$, and we attain a simpler hybrid model
by separating states that change during communications from those
that change during computations. 
}

\subsection{Combined Hybrid System} \label{ss:combined}
We are now ready to form the combined hybrid system for analysis. First, we combine all~$\eta^i$ values into a single variable~$\eta:= \col (\eta^1, \eta^2, \dots, \eta^N)$, which is in~$\R^{nN}$. We define the state of the combined hybrid system as~$\xi = (x, \eta, \tau) \in \mathcal{X}$, where~$\mathcal{X}:=\R^n \times \R^{nN} \times \T$ and~$\T := [0, \tau_{\max} ]$. To simplify notation, let the functions~$h_i: \R^n \to \R^{n_i}$ be given by~$h_i(\eta^i) = \nabla_{i} L(\eta^i)$ for all~$i \in [N]$. We collect these together into the function~$h: \R^{nN} \to \R^n$, given by
\begin{equation} \label{eq:hdef}
h(\eta):= \col \bigl(h_1(\eta^1), \dots, h_N(\eta^N)\bigr).
\end{equation}

This definition leads to the combined hybrid system~$\mathcal{H}=(C,f,D,G)$ with 
\begin{align} 
    \dot{\xi} & = \begin{bmatrix}
    -h(\eta) \\
    \boldsymbol{0}_{nN} \\
    -1
    \end{bmatrix} =: f(\xi) \label{eq:hfirst} 
\end{align}
for every~$\xi \in C:= \mathcal{X}$. 
Similar to the hybrid subsystems, when~$\tau= 0$, all agents undergo a jump. When a jump occurs,~$x$ remains constant,~$\eta^i$ is mapped with~$\eta^{i, +} = x$ for all~$i \in [N]$, and~$\tau^+ \in [\tau_{\min}, \tau_{\max}]$. Formally, for each~$\xi \in D:= \lbrace \xi \in \mathcal{X} : \tau = 0 \rbrace$, we define
the jump map~$G$ as 
 \begin{align}
    \xi^{+} &\in \begin{bmatrix}
    x \\
    \col (x, \dots, x ) \\
    [\tau_{\min}, \tau_{\max}]
    \end{bmatrix} =: G(\xi), \,
    \textnormal{dom } G = \R^n \times \R^{nN} \times \R. 
    \label{eq:hlast} 
 \end{align}

\new{
The connections between~\eqref{eq:hfirst}-\eqref{eq:hlast} and Problem~\ref{prob:first} are as follows. 
The block coordinate descent law that we use only has agent~$i$ perform computations
on~$x_i$ and~\eqref{eq:dynamics} shows this for agent~$i$. 
The first entry of~$f$ in~\eqref{eq:hfirst} is~$-h(\eta)$, and this entry replicates the dynamics in~\eqref{eq:dynamics}
for all~$i \in [N]$.
This entry of~$f$ therefore encodes that, for all~$i \in [N]$, agent~$i$ updates the value of~$x_i$ by moving
it along a negative gradient flow. 
Agent~$i$ stores a full copy of the decision variable~$x$ onboard 
itself in the state~$\eta^i$. 
Agent~$i$ does not update~$x_j$ for any~$j \neq i$
during its computations, and agent~$i$ therefore has~$\dot{\eta}^i_j = 0$ for all~$j \in [N] \backslash \{i\}$.
It also has~$\dot{\eta}^i_i = 0$ because~$\eta^i$ only stores values that have been communicated,
and the value of~$x_i$ is copied over to~$\eta^i_i$ only when~$x_i$ is communicated to other agents. 
These zero derivatives are encoded for all
agents in the second entry of~$f$ in~\eqref{eq:hfirst}, which is a zero vector of the appropriate size.
The third entry of~$f$ in~\eqref{eq:hfirst} is equal to~$-1$ to encode the fact that the timer~$\tau$
counts down with unit rate until it reaches zero. 
}

\new{
The jump map~$G$ in~\eqref{eq:hlast} is triggered when~$\tau = 0$ is reached. 
The first entry in~$G$ is~$x$, which shows that jumps do not change the decision variables that agents compute; 
the state~$x$ is simply set equal to its current value, which encodes no change. 
The second entry in~$G$ encodes the fact that agents communicate by 
broadcasting their states when~$\tau = 0$.
Mathematically, for all~$i \in [N]$, agent~$i$ sets~$\eta^i = x$ and
this step performs two distinct operations. First,
it sets~$\eta^i_i = x_i$, where~$x_i$ is the value
of agent~$i$'s decision variables when~$\tau = 0$ was reached.
This step ensures that this value of~$x_i$ will be used in
agent~$i$'s computations that will occur after the jump.
Second, this step sets~$\eta^i_j = x_j$ for~$j \in [N] \backslash \{i\}$,
which models agent~$j$ communicating to agent~$i$
the value of~$x_j$ from the time at which~$\tau = 0$ was reached. 
These new entries of~$\eta^i$ will also be used in agent~$i$'s updates after the jump is completed.
The third entry in~$G$ simply resets the value of~$\tau$ to some value
in the interval~$[\tau_{\min}, \tau_{\max}]$, and this reset deliberately allows
for some non-determinism in the value that~$\tau$ is reset to. 
}

\new{
This hybrid system models the information that drives all agents' computations and communications,
and our forthcoming analyses will derive the properties of this model, including convergence to minimizers. 
}

\subsection{Hybrid Basic Conditions} \label{ss:basicconditions}
We now demonstrate that~$\mathcal{H}$ meets the hybrid basic conditions and is well-posed.
\begin{lemma} \label{lem:hbasic} 
Let~$L$ satisfy Assumption~\ref{as:L}. Then the hybrid system given by~$\mathcal{H} = (C,f,D,G)$ and whose data is defined in~\eqref{eq:hfirst} and~\eqref{eq:hlast} satisfies the hybrid basic conditions from Definition~\ref{def:hybridbasic} and is nominally well-posed as a result.
\end{lemma}
\noindent \emph{Proof: }The sets~$C$ and~$D$ are closed subsets of~$\R^n \times \R^{nN} \times \R$ by definition. Due to our assumption that~$\nabla L$ is continuous,~$f$ is a continuous function 
from~$C$ to~$\R^n \times \R^{nN} \times \R$. By construction,~$G$ is outer semicontinuous and locally bounded relative to~$D$. Finally,~$D \subset \textnormal{dom } G$ because~$\textnormal{dom } G = \R^n \times \R^{nN} \times \R$ from~\eqref{eq:hlast}. 
\hfill $\qed$

\subsection{Existence of Solutions}
We denote solutions to~$\mathcal{H}$ by~$\phi$, which we parameterize by~$(t,j) \in \R_{\geq 0} \times \N$, where~$t$ denotes the ordinary (continuous) time, and~$j$ is a natural number that denotes the jump (discrete) time. Here,~$j$ is the cumulative number of jumps the agents have performed. Per Definition 2.3 in~\cite{goebel12},~$\textnormal{dom } \phi \subset \R_{\geq 0} \times \N$ is a~\emph{hybrid time domain} if for all~$(T,J) \in \textnormal{dom } \phi$, the set~$\textnormal{dom } \phi \cup \left([0,T]\times \{0,1,\dots,J\}\right)$ can be written as~$\bigcup_{j=0}^{J-1} ([t_j,t_{j+1}],j)$ for some finite sequence of times~$0 = t_0 \leq t_1 \leq \dots \leq t_J$. We say that a solution~$\phi$ is~\emph{complete} if~$\textnormal{dom } \phi$ is unbounded and we denote the set of all maximal solutions to~$\mathcal{H}$ as~$\mathcal{S}_\mathcal{H}$. A solution~$\phi$ to~$\mathcal{H}$ is called maximal if it cannot be extended further. 
In addition to being well-posed, there exists a nontrivial solution to~$\mathcal{H}$.
Below, in Proposition~\ref{prop:ges}, we will show that all
maximal solutions are complete. Toward doing so, we have
the following lemma. 

\begin{lemma}[Completeness of Solutions] \label{lem:solutions}
Let Assumption~\ref{as:L} hold and consider the hybrid system~$\mathcal{H}$ defined in~\eqref{eq:hfirst} and~\eqref{eq:hlast}. Let~$\tau_{\min}$ and~$\tau_{\max}$ be such that~$0 < \tau_{\min} \leq \tau_{\max}$. 
Then there exists a nontrivial solution to~$\mathcal{H}=(C,f,D,G)$ from every initial point in~$C \cup D$. 
Additionally, every maximal solution~$\phi$ to the hybrid system~$\mathcal{H}$ is non-Zeno and complete.
\end{lemma}
\noindent \emph{Proof: }See Appendix~\ref{app:hsm}. \hfill $\qed$

The next section also analyzes the stability of~$\mathcal{H}$, and as a preliminary result we have the following
lemma on how system trajectories evolve during flow intervals. 

\begin{lemma} \label{lem:contract}
Consider the hybrid system~$\mathcal{H}=(C,f,D,G)$ with data given 
in~\eqref{eq:hfirst} and~\eqref{eq:hlast}. Pick a solution~$\phi = (\phi_{x},\phi_{\eta},\phi_{\tau})$ to~$\mathcal{H}$. For each~$I^j := \{ t: (t,j) \in \textnormal{dom } \phi \}$ with nonempty interior and with~$t_{j+1} > t_j$ such that~$[t_j, t_{j+1}] = I^j$, we have
 \begin{align}
     \phi_{\eta^i}(t,j) &= \phi_{\eta^i}(t_j,j) \label{eq:z2update} \\
     \phi_{x_i}(t,j) &= \begin{cases} \phi_{\eta^i_i}(t_j,j) - (t-t_j) \nabla_i L\bigl(\phi_{\eta^i}(t_j,j)\bigr) &j \geq 1 \\
     \phi_{x_i}(0,0) - t \nabla_i L\bigl(\phi_{\eta^i}(0,0)\bigr) &j =0
     \end{cases} \label{eq:z1update}
 \end{align}
for all~$t \in (t_j, t_{j+1})$, where~$t_j$ denotes the continuous time at which the most recent jump~$j$ was performed.
\end{lemma}
\noindent \emph{Proof: } Given~$t\in(t_j, t_{j+1})$, the solution~$\phi = ( \phi_{x}, \phi_{\eta}, \phi_{\tau})$ has flowed some distance given by~$(t-t_j)\dot{\phi}$,  
where~$\dot{\phi}$ is constant due to the sample-and-hold
methodology.
Applying our definition of~$f$ in~\eqref{eq:hfirst} gives~\eqref{eq:z2update} and~\eqref{eq:z1update}.  \hfill $\qed$

\section{Stability Analysis} \label{sec:stability}
In this section, we define the convergence set~$\mathcal{A}$ and propose a Lyapunov function in Lemma~\ref{lem:comparison}. As an interim result, we show that if all agents initialize with the same state values, then their total distance from~$\mathcal{A}$ is monotonically decreasing for all objective functions that satisfy Assumptions~\ref{as:L} and~\ref{as:pl}. Next, Proposition~\ref{prop:ges} expands this interim result and bounds the distance from~$\mathcal{A}$ for all hybrid time~$(t,j)$ where~$j \geq 1$, regardless of agents' initialization. Finally,
Theorem~\ref{thm:3} removes the condition~$j \geq 1$
and establishes global exponential stability. 

\subsection{Convergence Set}
By Assumption~\ref{as:pl}, the set of minimizers~$\X^*$ 
is non-empty and may contain more than one element. Following from the properties of~$L$, namely Assumption~\ref{as:pl}, the algorithm has converged to a minimizer~$x^* \in \X^*$ of~$L$ if and only if it has reached a stationary point, i.e., a point at which~$\nabla L$ is zero. Given a complete solution~$\phi = (\phi_{x},\phi_{\eta},\phi_{\tau})$ to the hybrid system~$\mathcal{H}$, we seek to ensure that~$\lim_{t+j \to \infty} \nabla L\bigl(\phi_{x}(t,j)\bigr) = \boldsymbol{0}_{n}$ and~$\lim_{t+j \to \infty} \nabla L\bigl(\phi_{\eta^i}(t,j)\bigr) = \boldsymbol{0}_{n}$ for all~$i \in [N]$. This is equivalent to a set stability problem where the convergence set for the hybrid system~$\mathcal{H}$ is given by
\begin{align}
    \mathcal{A} &:= \lbrace \xi = (x, \eta, \tau) \in \mathcal{X} : \nabla L(x) = \boldsymbol{0}_n, \nabla L(\eta^i) = \boldsymbol{0}_n, \tau \in [0, \tau_{\max}],  \textnormal{for all } i \in [N] \rbrace \\
    & = \X^* \times \Bigl( \X^* \Bigr)^N \times[0, \tau_{\max}]. \label{eq:Adef}
\end{align}
Given a vector~$\xi=(x, \eta, \tau) \in \mathcal{X}$, let~$\hat{x}^0$ be the closest element of~$\X^*$ to~$x$, and let~$\hat{x}^i$ be the closest element of~$\X^*$ to~$\eta^i$ for each~$i \in [N]$. Formally, given~$\xi = (x, \eta, \tau)$, the points~$\hat{x}^0$ and~$\hat{x}^i$ are defined as
\begin{align}
    \hat{x}^0 := \argmin_{x^* \in \X^*}  \norm{x - x^*} \qquad \textnormal{ and } \qquad
    \hat{x}^i := \argmin_{x^* \in \X^*}  \norm{\eta^i - x^*} \textnormal{ for all } i \in [N]. 
\end{align}
Using these definitions, the squared distance from~$\xi$ to~$\mathcal{A}$ is given by~$\bigl| \xi\bigr|^2_\mathcal{A} := \norm{x - \hat{x}^0}^2 + \sum_{i \in [N]} \norm{\eta^i - \hat{x}^i}^2$. 

For all~$\xi = (x, \eta, \tau) \in \mathcal{A}$, the definition of~$\mathcal{A}$ does not immediately imply that~$x, \eta^1, \dots, \eta^N$ 
all converge to the same point~$x^* \in \X^*$. However, when combined with the dynamics of our hybrid system~$\mathcal{H}$, this
convergence property is guaranteed.

\begin{lemma}
Consider the hybrid system~$\mathcal{H}$ defined in~\eqref{eq:hfirst} and~\eqref{eq:hlast}. Let~$\mathcal{A}$ be as defined in~\eqref{eq:Adef}. For each maximal solution~$\phi$ to~$\mathcal{H}$, if~$\phi(t,j) \in \mathcal{A}$ for~$(t,j) \in \textnormal{dom }\phi$ with~$t \geq \tau_{\max}$, then 
    $\phi_{x}(t,j)  =\phi_{\eta^1}(t,j) = \dots = \phi_{\eta^N}(t,j) \in \X^*$.
\end{lemma}
\emph{Proof: } Following the definition of~$\mathcal{A}$, the condition~$\phi(t,j) \in \mathcal{A}$ implies both that $\nabla L\bigl(\phi_x (t,j)\bigr) = \boldsymbol{0}_n$ and~$\nabla L\bigl(\phi_{\eta^i}(t,j)\bigr) = \boldsymbol{0}_n$ for all~$i$. Note that because~$t \geq \tau_{\max}$, agents have performed at least one jump. Now, for the sake of contradiction, suppose that~$\phi_{\eta^i}(t,j) \neq \phi_x(t,j)$ for at least one~$i \in [N]$. Then there exists at least one entry~$\ell$ of~$\phi_x(t,j)$ such that~$\phi_{x_\ell}(t,j) \neq \phi_{\eta^i_\ell}(t,j)$. Because agents have performed at least one jump, it holds that~$\phi_{\eta^i_\ell}(t,j) = \phi_{\eta^\ell_\ell}(t,j)$. Combining this equality with~\eqref{eq:z2update} and~\eqref{eq:z1update} provides the relationship
    $\phi_{x_\ell}(t,j) = \phi_{\eta^i_\ell}(t,j) - (t-t_j) \nabla_\ell L\bigl(\phi_{\eta^\ell}(t,j)\bigr)$.
To satisfy the condition~$\phi_{x_\ell}(t,j) \neq \phi_{\eta^i_\ell}(t,j)$, it is necessary that~$\nabla_\ell L\bigl(\phi_{\eta^\ell}(t,j)\bigr) \neq 0$, 
which contradicts the hypothesis that~$\phi_{\eta^{\ell}}(t, j) \in \X^*$. Then~$\phi_{\eta^i}(t,j) = \phi_{x}(t, j) \in \X^*$ for all~$i \in [N]$ and~$t \geq \tau_{\max}$. 
\hfill $\qed$

\subsection{Bounds on the Lyapunov Function}
Central to proving our main result is a Lyapunov function that is bounded above and below by~$\mathcal{K}_\infty$ comparison functions~$\alpha_1$ and~$\alpha_2$, given next. 

\begin{lemma} \label{lem:comparison}
Let Assumptions~\ref{as:L} and~\ref{as:pl} hold. Let~$V : \mathcal{X} \to \R_{\geq 0}$ be a Lyapunov function candidate for the hybrid system~$\mathcal{H}= (C,f,D,G)$ defined 
in~\eqref{eq:hfirst} and~\eqref{eq:hlast}, given by
\begin{equation}
    V(\xi) = \bigl(L(x) - L^*\bigr) +  \sum_{i \in [N]} \bigl(L(\eta^i) - L^*\bigr),
\end{equation}
for all~$\xi=(x, \eta, \tau) \in \mathcal{X}$, where~$L$ is the objective function and~$L^*$ is from~\eqref{eq:lstar}. Then, there exist~$\alpha_1, \alpha_2 \in \mathcal{K}_\infty$ such that~$\alpha_1(| \xi|_\mathcal{A}) \leq V(\xi) \leq \alpha_2(| \xi|_\mathcal{A})$ for all~$\xi \in C \cup D \cup G(D)$. In particular, for all~$s\geq 0$,~$\alpha_1$ and~$\alpha_2$ are given by
\begin{align}
    \alpha_1(s) := \frac{\beta}{2} s^2 \quad \textnormal{and} \quad \alpha_2(s) := \frac{K}{2} s^2.
\end{align}
\end{lemma}
\noindent \emph{Proof: } See Appendix~\ref{app:stability}. \hfill $\qed$

\subsection{Global Exponential Stability}
We first bound the distance to a minimizer of~$L$ over time for a class of initial conditions in Proposition~\ref{prop:conv}. This result is then expanded to include all possible solutions and initial conditions in Proposition~\ref{prop:ges}, which characterizes the convergence
of trajectories after the first jump. Then, 
Theorem~\ref{thm:3} extends Proposition~\ref{prop:ges} to 
all times. 
We first consider the case where agents all agree at initialization and~$\phi_x(0,0) = \phi_{\eta^i}(0,0)$ for all~$i \in [N]$.
\begin{proposition} \label{prop:conv}
Let Assumptions~\ref{as:L} and~\ref{as:pl} hold and consider the hybrid system~$\mathcal{H}$ defined in~\eqref{eq:hfirst} and~\eqref{eq:hlast}. Let~$\mathcal{A}$ be as defined in~\eqref{eq:Adef} and let~$\tau_{\min}$ and~$\tau_{\max}$ be such that~$0 < \tau_{\min} \leq \tau_{\max} < \frac{1}{K}$, where~$K$ is the Lipschitz constant of~$\nabla L$
from Assumption~\ref{as:L}. Consider a maximal solution~$\phi$ to~$\mathcal{H}$ such that~$\phi_{x}(0,0) = \phi_{\eta^{i}}(0,0)$ for all~$i$ in~$[N]$. Then, for all~$(t,j) \in \textnormal{dom } \phi$, the following is satisfied:
\begin{equation}
    \bigl|\phi(t,j)\bigr|_\mathcal{A} \leq \sqrt{\frac{K}{\beta}} \exp \Big( - \frac{\beta}{N+1} ( 1 - K \tau_{\max}) t\Big) \bigl|\phi(0,0)\bigr|_\mathcal{A},
\end{equation}
where~$\beta$ is the PL constant of~$L$ from Assumption~\ref{as:pl} 
and~$1 - K \tau_{\max} > 0$ from the upper bound on~$\tau_{\max}$. 
\end{proposition}
\noindent \emph{Proof: } See Appendix~\ref{app:stability}. \hfill $\qed$

\new{
Proposition~\ref{prop:conv} shows that solutions to~$\mathcal{H}$ converge exponentially fast to a minimizer, 
and the exponent agrees (up to constants) with the convergence rate for centralized
continuous-time gradient descent on PL functions~\cite[Example 1]{wensing20}. 
It also agrees with the convergence rate for discrete-time multi-agent
minimization of PL functions~\cite[Theorem 1]{yazdani2021asynchronous} (again, up to constants). 
}

\newer{
Regarding the coefficient in Proposition~\ref{prop:conv}, the constant~$\sqrt{\frac{K}{\beta}}$ is tight in the sense that the bound
holds with equality at some times for some functions. For example,
consider~$L(x) = \frac{1}{2}\|x\|^2$, which satisfies
Assumption~\ref{as:L} with constant~$K=1$ and
satisfies Assumption~\ref{as:pl} with constant~$\beta=1$. 
For this choice of~$L$, 
at hybrid time~$(t,j) = (0, 0)$ the bound in Proposition~\ref{prop:conv} 
takes the form
\begin{equation}
\bigl|\phi(0,0)\bigr|_\mathcal{A} \leq \bigl|\phi(0,0)\bigr|_\mathcal{A},
\end{equation}
which shows that the constant~$\sqrt{\frac{K}{\beta}}$ cannot be made smaller. 
}

\newer{
The exponential term in Proposition~\ref{prop:conv}
contains the constant~$-\frac{\beta}{N+1}(1 - K\tau_{\max})$,
which differs from the constant~$-\beta$ that
appears in the convergence rate
of centralized continuous-time gradient descent~\cite[Example 1]{wensing}. 
If one used our methods to analyze a centralized setup, 
then the state~$\xi$
would only need to contain~$\eta^1$ and there
would be no need for~$\xi$ to contain~$x$ nor any~$\eta^i$ for~$i \geq 2$
because there would be no communications and no other agents to disagree with. 
Then one could use the Lyapunov
function~$V(\xi) = L(\eta^1) - L(x^*)$, and re-tracing
the steps of Proposition~\ref{prop:conv} with
this choice of~$V$ would eliminate
the factor of~$\frac{1}{N+1}$. 
Then the constant in the exponential term 
would take the form~$-\beta(1 - K\tau_{\max})$, and
the value of~$\tau_{\max}$ could be chosen arbitrarily
small to drive this constant arbitrarily close to~$-\beta$.
Through these changes, Proposition~\ref{prop:conv} 
comes arbitrarily close to recovering the convergence
rate for centralized continuous-time gradient descent
with PL functions. 
}

In practice, this preliminary result is useful when agreeing on initial values is easy to implement. However, we wish to account for all possible initialization scenarios. 
When agents disagree on initial conditions, their computations are not guaranteed to decrease the distance to minimizers
until after the first jump. The next result establishes an exponential upper bound from each initial condition,
including those at which agents disagree, and it accounts for the possible increase in the distance to the set of minimizers
that can occur before the first jump. 
It results in a larger bound than the one in Proposition~\ref{prop:conv}, as the following result 
presents. 

\begin{proposition}[Exponential bound for~$j \geq 1$] \label{prop:ges}
Let Assumptions~\ref{as:L} and~\ref{as:pl} hold and consider the hybrid system~$\mathcal{H}$ defined in~\eqref{eq:hfirst} and~\eqref{eq:hlast}. Let~$\mathcal{A}$ be as defined in~\eqref{eq:Adef} and choose~$\tau_{\min}$ and~$\tau_{\max}$ such that~$0 < \tau_{\min} \leq \tau_{\max} < \frac{1}{K}$, where~$K$ is the Lipschitz constant 
of~$\nabla L$ from Assumption~\ref{as:L}. 
For each maximal solution~$\phi$ and for all~$(t,j) \in \textnormal{dom } \phi$ such that~$j \geq 1$, the following is satisfied:
\begin{equation}
\bigl|\phi(t,j)\bigr|_\mathcal{A} \leq \sqrt{\frac{2K(N\!+\!1)}{\beta}} \exp \Big( - \frac{\beta}{N+1} ( 1 - K \tau_{\max}) t\Big) \bigl|\phi(0,0)\bigr|_\mathcal{A}, 
\end{equation}
where~$\beta$ is the PL constant of~$L$ from Assumption~\ref{as:pl} 
and~$1 - K \tau_{\max} > 0$ from the upper bound on~$\tau_{\max}$.
\end{proposition}
\noindent \emph{Proof: } See Appendix~\ref{app:stability}. \hfill $\qed$

Of course, for global exponential stability we must show
the exponential convergence to minimizers of all
trajectories from all initial conditions for all
times, not only for~$j \geq 1$. Accordingly, the 
following theorem does so and gives our main result
on global exponential stability. 

\begin{theorem} \label{thm:3}
Let Assumptions~\ref{as:L} and~\ref{as:pl} hold and consider the hybrid system~$\mathcal{H}$ defined in~\eqref{eq:hfirst} and~\eqref{eq:hlast}. Let~$\mathcal{A}$ be as defined in~\eqref{eq:Adef} and choose~$\tau_{\min}$ and~$\tau_{\max}$ such that~$0 < \tau_{\min} \leq \tau_{\max} < \frac{1}{K}$, where~$K$ is the Lipschitz constant of~$\nabla L$
from Assumption~\ref{as:L}. 
Then, the set $\cal A$ is globally exponentially stable for $\cal H$ defined 
in~\eqref{eq:hfirst}-\eqref{eq:hlast}, 
namely, for each maximal solution $\phi$ to $\cal H$, 
for all~$(t, j) \in \textnormal{dom } \phi$, 
we have
\begin{equation}
\bigl|\phi(t,j)\bigr|_\mathcal{A} \leq 
\max\left\{
\frac{\sqrt{2}}{\exp(-\rho\tau_{\max})}, 
\frac{\sqrt{1 + 2K^2\tau_{\max}^2}}{\exp(-\rho\tau_{\max})}, 
\sqrt{\frac{2K(N+1)}{\beta}}
\right\}\exp(-\rho t)\bigl|\phi(0,0)\bigr|, 
\end{equation}
where~$\beta$ is the PL constant of~$L$ from Assumption~\ref{as:pl}, we have
$1 - K \tau_{\max} > 0$ from the upper bound on~$\tau_{\max}$, and
$\rho = \frac{\beta}{N+1}(1 - K\tau_{\max})$. 
\end{theorem}
\noindent \emph{Proof: }
By Lemma~\ref{lem:solutions} we know that any
maximal solution~$\phi$ is also complete. 
Consider any~$t$ such that~$(t,0) \in \textnormal{dom } \phi$. We first seek to bound~$\bigl|\phi(t,0)\bigr|_\mathcal{A}$ with some constant. Define~$\bar{x}^0 := \argmin_{x^* \in \X^*}  \norm{\phi_{x} (t, 0) - x^*}$. Note that~$\phi_{\eta^i}(t, 0) = \phi_{\eta^i}(0, 0)$ for all~$i \in [N]$ and define~$\hat{x}^i := \argmin_{x^* \in \X^*}  \norm{\phi_{\eta^i}(t, 0) - x^*} = \argmin_{x^* \in \X^*}  \norm{\phi_{\eta^i}(0, 0) - x^*}$ for all~$i \in [N]$. We begin by expanding~$\bigl|\phi(t,0)\bigr|^2_\mathcal{A}$, where
\begin{align}
    \bigl|\phi(t,0)\bigr|^2_\mathcal{A} &= \norm{\phi_{x} (t, 0) - \bar{x}^0}^2 + \sum_{i \in [N]} \norm{\phi_{\eta^i}(t, 0) - \hat{x}^i}^2. \label{eq:distance_exp3}
\end{align}
We now define~$\hat{x}^0 := \argmin_{x^* \in \X^*} 
\norm{\phi_{x}(0,0) - x^*}$. Note that by definition of~$\bar{x}^0$,
~$\norm{\phi_{x} (t, 0) - \bar{x}^0}^2 \leq 
\norm{\phi_{x} (t, 0) - \hat{x}^0}^2$. Along with~$\phi_{\eta^i}(t, 0) = \phi_{\eta^i}(0, 0)$ for all~$i \in [N]$, this allows us to rewrite~\eqref{eq:distance_exp3} as
\begin{align}
\bigl|\phi(t,0)\bigr|^2_\mathcal{A} &\leq  
\norm{\phi_{x} (t, 0) - \hat{x}^0}^2 
+ \sum_{i \in [N]} \norm{\phi_{\eta^i}(0, 0) - \hat{x}^i}^2.\label{eq:distrewrite3}
\end{align}

We first upper bound~$\norm{\phi_{x_i} (t, 0) - \hat{x}^0_i}^2$ by applying Lemma~\ref{lem:contract} and 
using~$\norm{a-b}^2 \leq 2\norm{a}^2 + 2\norm{b}^2$, resulting in
\begin{align}
   \norm{\phi_{x_i} (t, 0) - \hat{x}^{0}_i}^2 &= 
   \norm{\phi_{x_i} (0,0) - t \nabla_i L\bigl(\phi_{\eta^i} (0,0)\bigr) - \hat{x}^{0}_i}^2 \\
    &\leq 2 \norm{\phi_{x_i} (0,0) -  \hat{x}^{0}_i}^2 
    + 2 t^2 \norm{\nabla_i L\bigl(\phi_{\eta^i}  (0,0)\bigr) - \nabla_i L\bigl(\hat{x}^i \bigr)}^2 \\
    &\leq 2 \norm{\phi_{x_i} (0,0) -  \hat{x}^{0}_i}^2 + 2 K^2 \tau_{\max}^2 \norm{\phi_{\eta^i}(0,0) - \hat{x}^i}^2, \label{eq:xupperPL3}
\end{align}
where the first equality applies Lemma~\ref{lem:contract}, the first inequality uses~$\nabla L(\hat{x}^i) = 0$, and the final inequality applies Lemma~\ref{lem:lipschitz}.
Summing over all~$i$ on both sides of~\eqref{eq:xupperPL3} gives
\begin{align}
    \norm{\phi_{x} (t, 0) - \hat{x}^{0}}^2 
    \leq 2 \norm{\phi_{x} (0,0) -  \hat{x}^{0}}^2 + 2 K^2 \tau_{\max}^2 \sum_{i \in [N]} \norm{\phi_{\eta^i} (0,0) -  \hat{x}^i}^2. \label{eq:xupperPL4}
\end{align}
Applying this inequality to~\eqref{eq:distrewrite3} gives 
\begin{align}
    \bigl|\phi(t,0)\bigr|^2_\mathcal{A} &\leq 2 
    \norm{\phi_{x} (0,0) -  \hat{x}^{0}}^2 + \left( 1+ 2 K^2 \tau_{\max}^2 \right) 
    \sum_{i \in [N]} \norm{\phi_{\eta^i} (0,0) -  \hat{x}^i}^2 \\
    &\leq \max \{2, 1+ 2 K^2 \tau_{\max}^2 \} \bigl|\phi(0,0)\bigr|^2_\mathcal{A}.
\end{align}
Taking the square root and combining with Proposition~\ref{prop:ges} gives the final result.

\section{Robustness to Timing Errors} \label{sec:robustness}
In this section, we show that the hybrid system~$\mathcal{H}$ is robust to a class of model errors, in particular that it is robust to errors in the dynamics of the timer. 
By ``robust'' we mean that there exists a maximum nonzero perturbation level such that all solutions under such perturbations converge to a neighborhood of the set $\cal A$, where the size of the neighborhood depends on the size of the perturbation.
Formally, we consider perturbed timer dynamics of the form
\begin{equation} \label{eq:pdefs}
\dot{\tau}_p = -1 + \kappa, \quad \tau_p^+ \in [\tau_{\min} + \theta_{\min}, \tau_{\max} + \theta_{\max}],
\end{equation}
where~$\tau_p$ denotes the perturbed timer,~$\kappa \in (-\infty, 1)$ is a constant that models skew on the timer dynamics, and the terms~$\theta_{\min} \in \R$ and~$\theta_{\max} \in \R$ are perturbations to~$\tau_{\min}$ and~$\tau_{\max}$, respectively, that satisfy
\begin{equation} \label{eq:thetadefs}
0 < \tau_{\min} + \theta_{\min} \leq \tau_{\max} + \theta_{\max}.
\end{equation}
The full perturbation of the hybrid system model 
in~\eqref{eq:hfirst}-\eqref{eq:hlast} 
 has 
state vector 
denoted~$\xi_p = (x, \eta, \tau_p)$ 
whose flow dynamics are given by
\begin{equation}
\dot{\xi}_p = \left[\begin{array}{c} -h(\eta) \\ \boldsymbol{0}_{nN} \\
-1 + \kappa
\end{array}\right] 
=: f_p(\xi_p), 
\quad \quad \xi_p \in C_p := 
\R^n \times \R^{nN} \times [0, \tau_{\max} + \theta_{\max}],
\end{equation}
where~$h$ and~$\eta$ are from~\eqref{eq:hfirst}, and~$\kappa$ and~$\theta_{\max}$
are from~\eqref{eq:pdefs}. Its jump dynamics are given
by
\begin{equation}
\xi_p^{+} = \left[\begin{array}{c}
x \\
\textnormal{col}(x, \ldots, x) \\
\big[\tau_{\min} + \theta_{\min}, \tau_{\max} + \theta_{\max}\big]
\end{array}\right]
=: G_p(\xi_p), \quad
\quad \xi_p \in D_p := \{\xi_p \in \mathcal{X}_p : \tau_p = 0\}. 
\end{equation} 
We use
\begin{equation} \label{eq:hp}
\mathcal{H}_p := (C_p, f_p, D_p, G_p)
\end{equation}
to denote the 
full hybrid system with the perturbed timer dynamics. 

To enable the desired robustness
result, we first require the following assumption.

\begin{assumption} \label{as:Xcompact}
The set of optimizers~$\mathcal{X}^*$ is compact.
\end{assumption}

Assumption~\ref{as:Xcompact} is required here so that agents' computations do not converge to a solution that is arbitrarily far away from their current iterates. It is known to hold under mild conditions, such as the condition that the objective function~$L$ 
is coercive~\cite[Proposition 2.1.1]{bertsekas03}.
We have the following robustness result.

\begin{theorem}
Let Assumptions~\ref{as:L}-\ref{as:Xcompact} hold. 
Then, for the hybrid system~$\mathcal{H}_p$ from~\eqref{eq:hp}, 
there exists a function~$\beta \in \mathcal{KL}$ such that for
every~$\epsilon > 0$ there exists~$\rho^* \in (0, \infty)$ such that
if~$\max\big\{|\kappa|, |\theta_{\min}|, |\theta_{\max}|\big\} \leq \rho^*$, 
then each solution~$\phi_p$
to~$\mathcal{H}_p$ satisfies
\begin{equation}
|\phi_p(t, j)|_{\mathcal{A}} \leq
\beta\big(|\phi_p(0, 0)|_{\mathcal{A}}, t+j\big) + \epsilon
\end{equation}
for all~$(t, j) \in \textnormal{dom} \,\phi$, where~$\mathcal{A}$ is given
in~\eqref{eq:Adef}. 
\end{theorem}
\emph{Proof:} 
By Lemma~\ref{lem:hbasic},
the nominal hybrid
system~$\mathcal{H}$ satisfies
the hybrid basic conditions. 
By Assumption~\ref{as:Xcompact}
the set~$\mathcal{X}^*$ is 
compact and thus the set~$\mathcal{A}$ is compact. 
The
function~$x \mapsto |x|_{\mathcal{A}}$
is a proper indicator
for the set~$\mathcal{A}$
viewed as a subset of~$\R^n$.
That is, we have
both~$|x|_{\mathcal{A}} \to \infty$
as~$|x| \to \infty$
and~$|x|_{\mathcal{A}} = 0$
if and only if~$x \in \mathcal{A}$.
Theorem~\ref{thm:3} provides a
class-$\mathcal{KL}$ function~$\beta$ such that
\begin{equation}
|\phi(t,j)|_{\mathcal{A}} \leq \beta\big(|\phi(0,0)|_{\mathcal{A}}, t+j\big)
\end{equation}
for each solution~$\phi$
to~$\mathcal{H}$
and all~$(t, j)$
in the domain of~$\phi$.

The system~$\mathcal{H}_p$ in~\eqref{eq:hp} can be modeled as a~$\rho$-perturbation 
of~$\mathcal{H}$ in~\eqref{eq:hfirst}-\eqref{eq:hlast} — see Section 2.3.5 
of~\cite{sanfelice21} and~\cite[Exercise 25]{sanfelice21}. 
The perturbed flow map is equal to the nominal one plus $(0,0,\kappa)$. Then, given $\kappa \in (-\infty, 1)$, there exists $\rho^a >0$ such that $f_p(\xi) \subset f(\xi) + (0,0,\rho^a)\mathbb{B}$ for all $\xi$, where~$\mathbb{B}$ denotes the closed 
Euclidean unit ball.
Given $\theta_{\min}$ and $\theta_{\max}$ satisfying~\eqref{eq:thetadefs} there exists $\rho^b > 0$ such that $C_p \subset C+\rho^b \mathbb{B}$ and $[\tau_{\min}+\theta_{\min},\tau_{\max}+\theta_{\max}] \subset [\tau_{\min},\tau_{\max}] + \rho^b\mathbb{B}$.  
The jump map~$G_p$ satisfies~$G_p(\xi) \subset G(\xi) + (0, 0, \rho^b)\mathbb{B}$
for all~$\xi$, and the jump set is simply~$D_p = D$. 
The perturbed quantities~$f_p$, $C_p$, $G_p$, and~$D_p$
satisfy Assumption 3.25 in~\cite{sanfelice21} by inspection. 
Then all conditions in~\cite[Theorem 3.26]{sanfelice21} 
hold with $\rho := \max\{\rho^a,\rho^b\}$, and the theorem follows. \hfill $\qed$

\section{Numerical Validation} \label{sec:sim}
\begin{figure}[ht!]
\centering
\includegraphics[width=5.8cm]{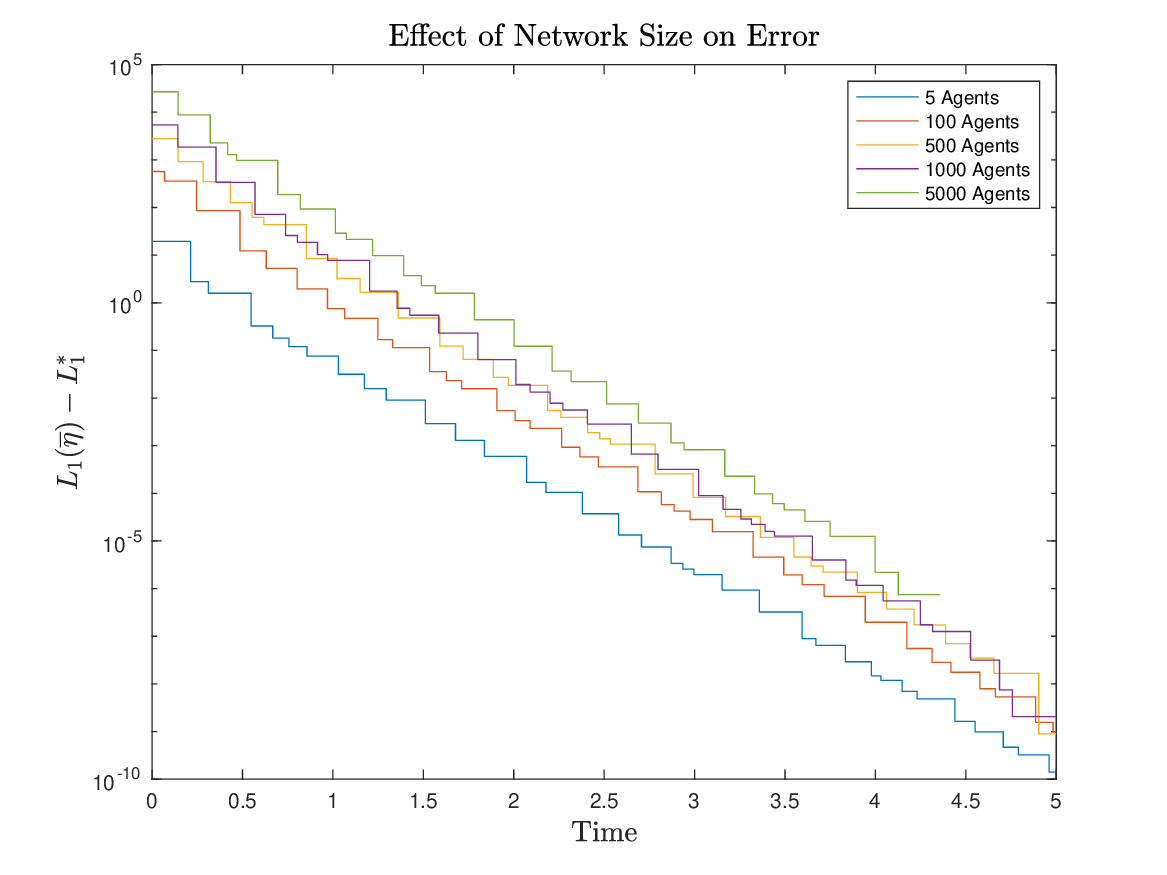}
\caption{Effect of network size on convergence for a strongly convex quadratic program. This is a semi-log plot, so 
straight lines imply exponential convergence. The horizontal axis is continuous time, and jumps are demonstrated by the sudden drops as they occur in discrete time. We see that exponential convergence is attained for all network sizes, demonstrating the scalability of our algorithm.}
\label{fig:qp_networksize}
\end{figure}

\begin{figure}[ht!]
\centering
     \begin{subfigure}[b]{0.48\textwidth}
         \centering
         \includegraphics[width=\textwidth]{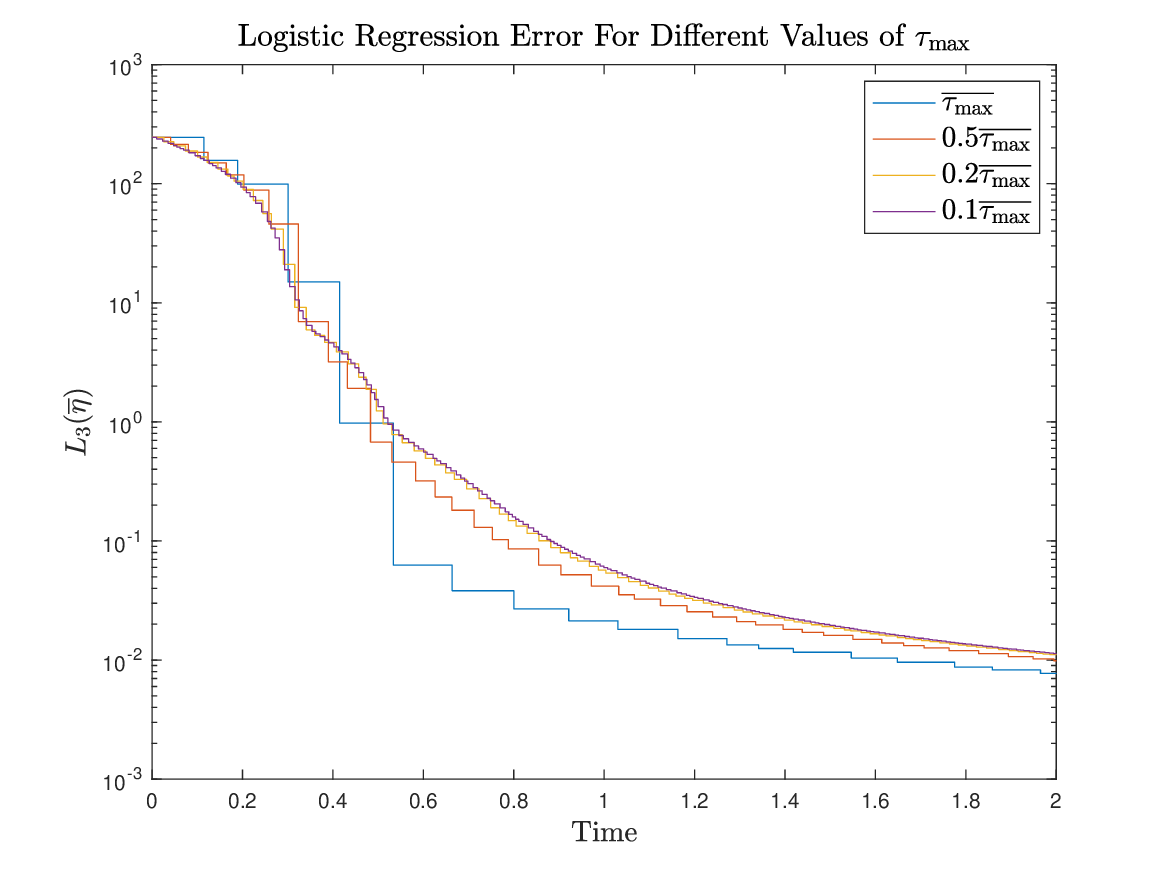}
         \caption{Effect of~$\tau_{\max}$ on convergence demonstrating that larger values of~$\tau_{\max}$ may help convergence. The largest value of~$\tau_{\max}$, indicated with the blue line, achieves the best performance.}
         \label{fig:lr_taumax}
     \end{subfigure} 
     \hfill
    \begin{subfigure}[b]{0.48\textwidth}
         \centering
         \includegraphics[width=\textwidth]{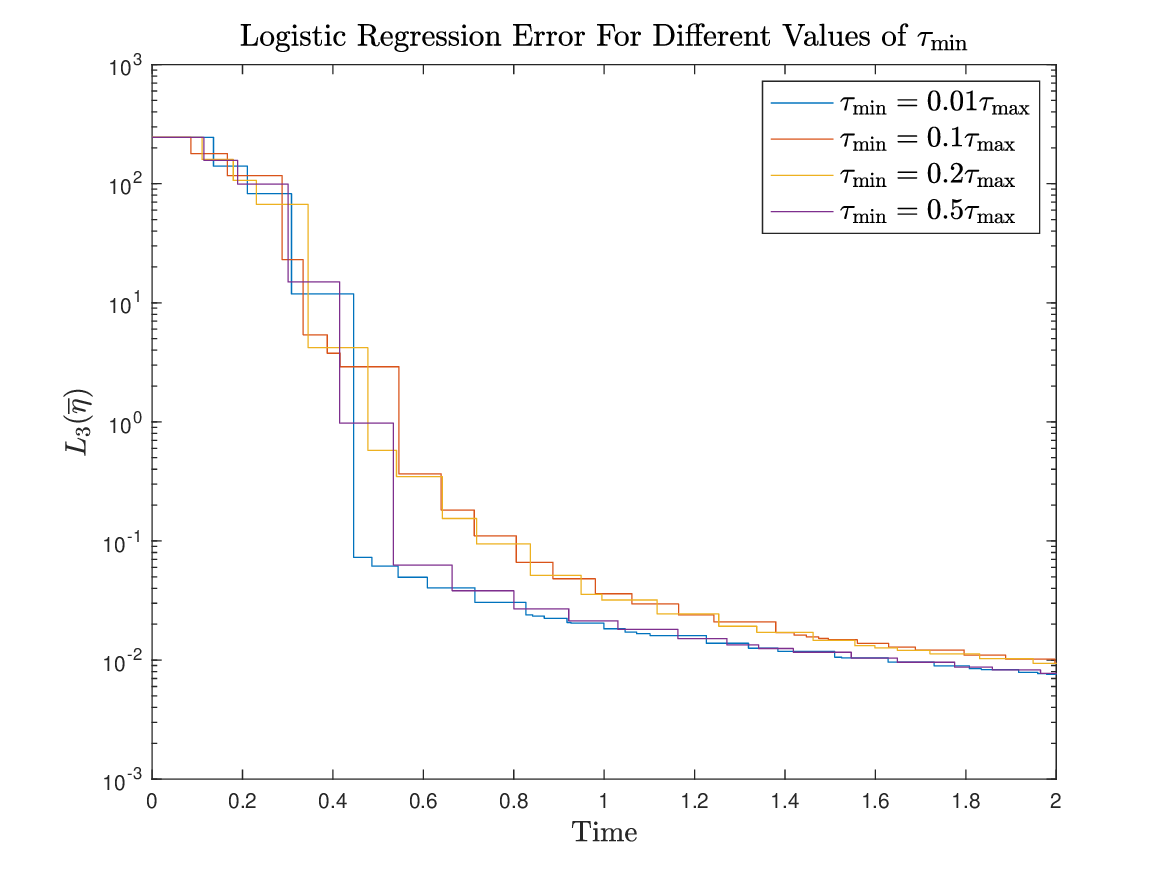}
         \caption{Effect of~$\tau_{\min}$ on convergence demonstrating the insignificant impact it appears to have on convergence. In contrast to the choice of~$\tau_{\max}$, for at least some problems, the choice of~$\tau_{\min}$ does not seem significant. The plot here shows that both the smallest (blue line) and largest (purple line) choices for~$\tau_{\min}$ have similar performance.}
         \label{fig:lr_taumin}
     \end{subfigure}
     \hfill
\caption{Effects of~$\tau_{\max}$ and~$\tau_{\min}$ on convergence for the logistic regression problem in Application 3.}
\label{fig:taumaxtaumin}
\end{figure}

\begin{figure}[ht!]
\centering
\includegraphics[width=5.8cm]{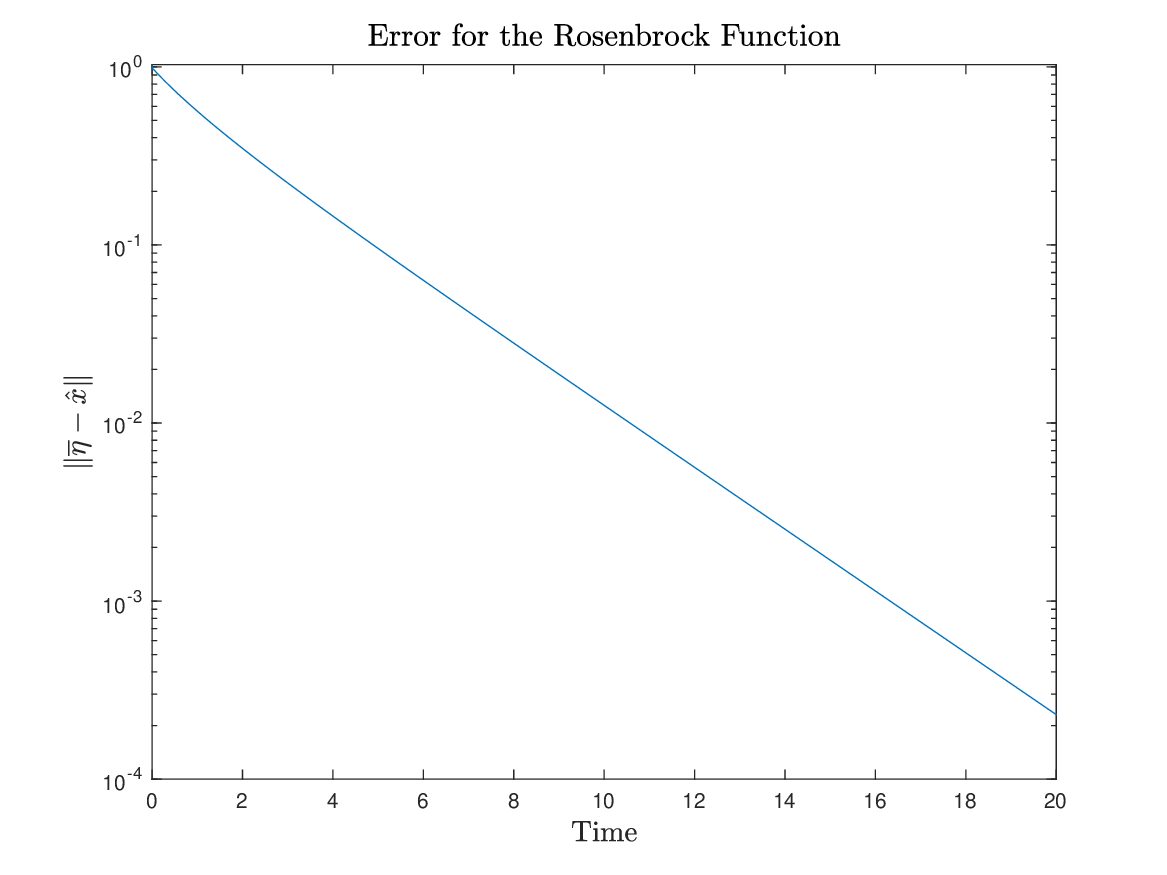}
\caption{Distance from the optimizer for the Rosenbrock problem, which decreases exponentially. 
}
\label{fig:rose}
\end{figure}

Three different applications are considered in this section: quadratic programs, linear neural networks (inspired by~\cite{charles18}), logistic regression, and the Rosenbrock problem~\cite{rosenbrock60}. In all cases, the HyEq Toolbox (Version 2.04)~\cite{HyEqToolbox} was used for simulation. 
Code for all simulations is available on
GitHub\footnote{Figures~\ref{fig:qp_networksize}, \ref{fig:taumaxtaumin}, and~\ref{fig:rose} 
are from simulations generated with code 
at \url{http://github.com/kathendrickson/DistrHybridGD}. 
Figure~\ref{fig:newapplication1} is from simulations generated with code
at \url{https://github.com/corelabgt/DistrHybridGDComparison}. 
}.

\begin{application}[Quadratic Program] \label{app:qp}
We consider~$N$ agents for the values $N \in \{5, 100, 500, 1000, 5000\}$. Each agent updates a scalar and they collaboratively minimize a quadratic function of the form
\begin{align} \label{eq:l1}
    L_1(x) := \frac{1}{2}x^{\top} Q x + b^{\top}x,
\end{align}
where~$x\in \R^N$,~$Q$ is an~$N \times N$ symmetric, positive definite matrix, and~$b\in [1,5]^N$. For all experiments,~$\tau_{\max} = \frac{1}{K + 0.001}$ and~$\tau_{\min} = \frac{1}{5} \tau_{\max}$.
\end{application}
The function~$L_1$ in~\eqref{eq:l1} from Application~\ref{app:qp} is strongly convex and smooth, hence Assumptions~\ref{as:L} and~\ref{as:pl} hold with parameters~$\beta = \min \lambda (Q)$ and~$K=\max \lambda (Q)$, where~$\lambda(Q)$ denotes the set of eigenvalues of the matrix~$Q$. 
 
We consider~$\beta = 2$ and~$K = 4$ and initialize~$\phi_{x} = \phi_{\eta^i}$ for all~$i\in [N]$. Matlab's Optimization Toolbox was used to find~$L_1^*$, which is compared to~$L_1 (\overline{\eta})$, the objective function evaluated at the shared value of~$\eta$, throughout the experiment. As shown in Figure~\ref{fig:qp_networksize}, expanding the network size does not have a significant impact on convergence. This demonstrates our algorithm's scalability and convergence that holds regardless of network size.

\new{
The preceding simulations compare problems
of variable size because~$x \in \mathbb{R}^N$
for~$N \in \{5, 100, 500, 1000, 5000\}$.
To assess scalability, 
we also compare problems of the same size
parallelized among different numbers of agents. 
We consider~$x \in \mathbb{R}^{5000}$ and
again consider~$N$ agents
for~$N \in \{5, 100, 500, 1000, 5000\}$,
where each agent is responsible for 
updating~$5000/N$ decision variables. The problem
setup and choices of problem parameters
are the same as in the preceding simulation runs, 
and the results are shown in Figure~\ref{fig:newapplication1}.
The curves for each run overlap almost entirely, and
there is negligible impact in changing the number
of agents across which the problem is partitioned, which
further demonstrates that the convergence rate of our algorithm 
scales well. 
\begin{figure}[ht!]
\centering
\includegraphics[width=5.8cm]{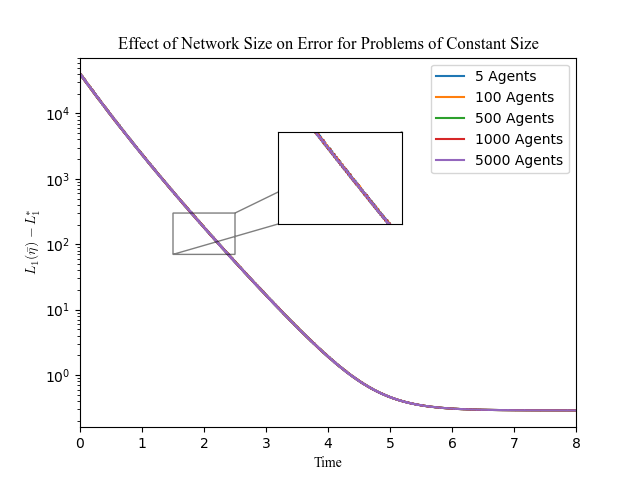}
\caption{
\new{
Additional runs of Application~1 in which the problem size
is fixed and the number of agents is variable.
Changes in the number of agents have a negligible impact
on the convergence rate, which further illustrates the
scalability of our algorithm. 
}
}
\label{fig:newapplication1}
\end{figure}
}

\begin{application}[Logistic Regression] \label{app:lr}
We consider~$N=100$ agents collaboratively minimizing a logistic regression cost function of the form
\begin{align}
    L_3(x) := \frac{1}{5} \sum_{i=1}^5 \log (1 + \exp{b_i a_i^{\top} x}),
\end{align}
where~$x \in \R^{N}$,~$b_i \in [0,10]$, and~$a_i \in \{0,1\}^{N}$ for~$i=1, \dots, 5$. The parameters~$\tau_{\max}$ and~$\tau_{\min}$ take various values that are shown on the plots below. 
\end{application}
While the logistic regression problem given by~$L_3$ is smooth (satisfying Assumption~\ref{as:L}) and convex, it is not strongly convex. 
However, according to~\cite[Section 2.3]{karimi2016linear},~$L_3$ satisfies the PL condition, which is Assumption~\ref{as:pl}, over any compact set. We therefore define the
set
$\big\{\xi = (x, \eta, \tau) : |\xi|_{\mathcal{A}} \leq |\phi(0, 0)|_{\mathcal{A}}\big\}$,
which is compact by construction, and we know that~$L_3$ satisfies Assumption~\ref{as:L} over this
set. Then, Theorem~\ref{thm:3} implies that 
system trajectories
remain in this set for all time, and thus our convergence reuslts apply to it.

A benchmark value for~$\tau_{\max}$ is set to~$\overline{\tau_{\max}} = \frac{1}{K + 0.001}$. As shown in Figure~\ref{fig:lr_taumax}, smaller values of~$\tau_{\max}$ may result in slower convergence. While it may seem counter-intuitive that larger bounds on delays may~\emph{help} convergence, this observation echoes work in discrete-time optimization by two of the authors~\cite{hendrickson2021totally} 
that found it necessary to balance between 
(i) delaying communications to allow agents to make progress with their current state values and
(ii) communicating more often to reduce disagreements among agents. 
In this problem, the interpretation of this idea is that a larger~$\tau_{\max}$ allows agents'
computations to 
make progress toward a minimizer before their next communication, which produces a net benefit to the overall
convergence of the algorithm. 
Next, various choices of~$\tau_{\min}$ are compared. Results shown in Figure~\ref{fig:lr_taumin} suggest that the choice of~$\tau_{\min}$ does not significantly impact convergence.

\begin{application}[Rosenbrock Problem~\cite{rosenbrock60}] \label{appl:rose}
Consider~$N=2$ agents minimizing the Rosenbrock function given by
\begin{align}
L_4(x) := (1-x_1)^2 + 100(x_2 - x_1^2)^2,
\end{align}
where~$x = (x_1, x_2) \in \R^2$. 
Here, the values~$\tau_{\max} = 0.001$ and~$\tau_{\min} = \frac{1}{5} \tau_{\max}$ were used.
\end{application}
The non-convex problem given by~$L_4$ is often used as a benchmark problem for optimization algorithms as it is difficult to solve due to the problem's geometry. However, it satisfies Assumptions~\ref{as:L} and~\ref{as:pl} in the region~$[-1,1]^2$ and has a global optimum at the point~$(1,1)$~\cite{budhraja21}. 
Figure~\ref{fig:rose} plots the distance of the shared value~$\overline{\eta}$ from the global optimum at~$(1,1)$ throughout the simulation. While convergence is slower than the previous examples and the bound on~$\tau_{\max}$ is relatively small compared to previous examples, exponential convergence is still achieved.

\section{Conclusion} \label{sec:concl}
This paper presented a hybrid systems framework for analyzing continuous-time multi-agent optimization with discrete-time communications. Using this framework, we were able to establish that every maximal solution is complete, as well as the global exponential convergence of a block coordinate descent law to a minimizer of a smooth, possibly nonconvex, objective function that satisfies the PL inequality. Finally, three applications were considered with simulation results demonstrating the scalability and performance of our framework. 

\begingroup
\let\itshape\upshape
\bibliographystyle{plainnat}
\bibliography{sources}

@book{sra12,
  title={Optimization for Machine Learning},
  author={Sra, Suvrit and Nowozin, Sebastian and Wright, Stephen J},
  year={2012},
  publisher={MIT Press}
}

@book{bertsekas89,
 author = {Bertsekas, Dimitri P. and Tsitsiklis, John N.},
 title = {Parallel and Distributed Computation: Numerical Methods},
 year = {1989},
 publisher = {Prentice-Hall, Inc.},
 address = {USA}
}

@book{goebel12,
  title={Hybrid Dynamical Systems: Modeling, Stability, and Robustness},
  author={Goebel, R. and Sanfelice, R.G. and Teel, A.R.},
  year={2012},
  publisher={Princeton University Press},
  address={Princeton, NJ, USA}
}

@article{wensing,
    author = {Wensing, Patrick M. AND Slotine, Jean-Jacques},
    journal = {PLOS ONE},
    publisher = {Public Library of Science},
    title = {Beyond convexity—Contraction and global convergence of gradient descent},
    year = {2020},
    month = {08},
    volume = {15},
    url = {https://doi.org/10.1371/journal.pone.0236661},
    pages = {1-29},
    number = {8},
    doi = {10.1371/journal.pone.0236661}
}

@article{chai17,
  title={Analysis and design of event-triggered control algorithms using hybrid systems tools},
  author={Chai, Jun and Casau, Pedro and Sanfelice, Ricardo G},
  journal={International Journal of Robust and Nonlinear Control},
  volume={30},
  number={15},
  pages={5936--5965},
  year={2020},
  note={\url{https://doi.org/10.1002/rnc.5141}}
}

@inproceedings{HyEqToolbox, 
  author = {Sanfelice, Ricardo G. and Copp, David and Nanez, Pablo}, 
  title = {A Toolbox for Simulation of Hybrid Systems in Matlab/Simulink: Hybrid Equations {(HyEQ)} Toolbox}, 
  year = {2013}, 
  note = {\url{https://doi.org/10.1145/2461328.2461346}}, 
  pages = {101–106},
  booktitle = {Proceedings of the 16th International Conference on Hybrid Systems: Computation and Control} 
}

@article{su16,
  author  = {Weijie Su and Stephen Boyd and Emmanuel J. Cand{{\`e}}s},
  title   = {A Differential Equation for Modeling {N}esterov's Accelerated Gradient Method: Theory and Insights},
  journal = {Journal of Machine Learning Research},
  year    = {2016},
  volume  = {17},
  number  = {153},
  pages   = {1-43},  
}

@article{KIA2015,
title = {Distributed convex optimization via continuous-time coordination algorithms with discrete-time communication},
journal = {Automatica},
volume = {55},
pages = {254-264},
year = {2015},
note = {\url{https://doi.org/10.1016/j.automatica.2015.03.001}},
author = {Solmaz S. Kia and Jorge Cortés and Sonia Martínez},
}

@ARTICLE{rahili2017,
  author={S. {Rahili} and W. {Ren}},
  journal={IEEE Transactions on Automatic Control}, 
  title={Distributed Continuous-Time Convex Optimization With Time-Varying Cost Functions}, 
  year={2017},
  volume={62},
  number={4},
  pages={1590-1605},
  note={\url{https://doi.org/10.1109/TAC.2016.2593899}}
}

@ARTICLE{gharesifard2014,
  author={B. {Gharesifard} and J. {Cortés}},
  journal={IEEE Transactions on Automatic Control}, 
  title={Distributed Continuous-Time Convex Optimization on Weight-Balanced Digraphs}, 
  year={2014},
  volume={59},
  number={3},
  pages={781-786},
  note={\url{https://doi.org/10.1109/TAC.2013.2278132}}
}

@ARTICLE{lu2012,
  author={J. {Lu} and C. Y. {Tang}},
  journal={IEEE Transactions on Automatic Control}, 
  title={Zero-Gradient-Sum Algorithms for Distributed Convex Optimization: The Continuous-Time Case}, 
  year={2012},
  volume={57},
  number={9},
  pages={2348-2354},
  note={\url{https://doi.org/10.1109/TAC.2012.2184199}}
}

@article{PHILLIPS2019,
title = {Robust distributed synchronization of networked linear systems with intermittent information},
journal = {Automatica},
volume = {105},
pages = {323-333},
year = {2019},
note = {\url{https://doi.org/10.1016/j.automatica.2019.03.020}},
author = {Sean Phillips and Ricardo G. Sanfelice},
}

@INPROCEEDINGS{phillips2018,  
author={Phillips, Sean and Erwin, R. Scott and Sanfelice, Ricardo G.},  
booktitle={2018 Annual American Control Conference (ACC)},   
title={Robust Exponential Stability of an Intermittent Transmission State Estimation Protocol},   
year={2018},  
volume={},  
number={},  
pages={622-627},  
note = {\url{https://doi.org/10.23919/ACC.2018.8431144}}
}

@article{garg20,
  title={Fixed-time stable gradient flows: Applications to continuous-time optimization},
  author={Garg, Kunal and Panagou, Dimitra},
  journal={IEEE Transactions on Automatic Control},
  volume={66},
  number={5},
  pages={2002--2015},
  year={2020},
  note={\url{https://doi.org/10.1109/TAC.2020.3001436}}
}

@article{luo06,
  title={An introduction to convex optimization for communications and signal processing},
  author={Luo, Zhi-Quan and Yu, Wei},
  journal={IEEE Journal on Selected Areas in Communications},
  volume={24},
  number={8},
  pages={1426--1438},
  year={2006},
  note={\url{https://doi.org/10.1109/JSAC.2006.879347}}
}

@article{verscheure09,
  title={Time-optimal path tracking for robots: A convex optimization approach},
  author={Verscheure, Diederik and Demeulenaere, Bram and Swevers, Jan and De Schutter, Joris and Diehl, Moritz},
  journal={IEEE Transactions on Automatic Control},
  volume={54},
  number={10},
  pages={2318--2327},
  year={2009},
  note={\url{https://doi.org/10.1109/TAC.2009.2028959}}
}

@article{olfati07,
  title={Consensus and cooperation in networked multi-agent systems},
  author={Olfati-Saber, Reza and Fax, J Alex and Murray, Richard M},
  journal={Proceedings of the IEEE},
  volume={95},
  number={1},
  pages={215--233},
  year={2007},
  note={\url{https://doi.org/10.1109/JPROC.2006.887293}}
}

@article{cortes04,
  title={Coverage control for mobile sensing networks},
  author={Cortes, Jorge and Martinez, Sonia and Karatas, Timur and Bullo, Francesco},
  journal={IEEE Transactions on Robotics and Automation},
  volume={20},
  number={2},
  pages={243--255},
  year={2004},
  note={\url{https://doi.org/10.1109/TRA.2004.824698}}
}

@inproceedings{karimi2016linear,
  title={Linear convergence of gradient and proximal-gradient methods under the {P}olyak-{\L}ojasiewicz condition},
  author={Karimi, Hamed and Nutini, Julie and Schmidt, Mark},
  booktitle={Joint European Conference on Machine Learning and Knowledge Discovery in Databases},
  pages={795--811},
  year={2016},
  organization={Springer},
  note={\url{https://doi.org/10.1007/978-3-319-46128-1_50}}
}

@inproceedings{charles18,
  title={Stability and generalization of learning algorithms that converge to global optima},
  author={Charles, Zachary and Papailiopoulos, Dimitris},
  booktitle={International Conference on Machine Learning},
  pages={745--754},
  year={2018},
  organization={PMLR},
}

@article{rosenbrock60,
    author = {Rosenbrock, H. H.},
    title = {An Automatic Method for Finding the Greatest or Least Value of a Function},
    journal = {The Computer Journal},
    volume = {3},
    number = {3},
    pages = {175-184},
    year = {1960},
    note = {\url{https://doi.org/10.1093/comjnl/3.3.175}}
}

@inproceedings{budhraja21,
  title={Breaking the convergence barrier: Optimization via fixed-time convergent flows},
  author={Budhraja, Param and Baranwal, Mayank and Garg, Kunal and Hota, Ashish},
  booktitle={Proceedings of the AAAI Conference on Artificial Intelligence},
  volume={36},
  pages={6115--6122},
  year={2022},
  note={\url{https://doi.org/10.1609/aaai.v36i6.20559}}
}

@inproceedings{du2019gradient,
  title={Gradient descent finds global minima of deep neural networks},
  author={Du, Simon and Lee, Jason and Li, Haochuan and Wang, Liwei and Zhai, Xiyu},
  booktitle={International Conference on Machine Learning},
  pages={1675--1685},
  year={2019},
  organization={PMLR}
}

@inproceedings{fazel2018global,
  title={Global convergence of policy gradient methods for the linear quadratic regulator},
  author={Fazel, Maryam and Ge, Rong and Kakade, Sham and Mesbahi, Mehran},
  booktitle={International Conference on Machine Learning},
  pages={1467--1476},
  year={2018},
  organization={PMLR}
}

@article{lian2015asynchronous,
  title={Asynchronous parallel stochastic gradient for nonconvex optimization},
  author={Lian, Xiangru and Huang, Yijun and Li, Yuncheng and Liu, Ji},
  journal={Advances in Neural Information Processing Systems},
  volume={28},
  pages={2737--2745},
  year={2015}
}

@incollection{Chiang2009,
author="Chiang, Mung",
editor="Gao, David Y.
and Sherali, Hanif D.",
title="Nonconvex Optimization for Communication Networks",
bookTitle="Advances in Applied Mathematics and Global Optimization",
year="2009",
pages="137--196",
note={\url{https://doi.org/10.1007/978-0-387-75714-8_5}},
publisher={Springer},
}

@article{polyak1963gradient,
  title={Gradient methods for minimizing functionals},
  author={Polyak, Boris Teodorovich},
  journal={Zhurnal vychislitel'noi matematiki i matematicheskoi fiziki},
  volume={3},
  number={4},
  pages={643--653},
  year={1963},
  publisher={Russian Academy of Sciences, Branch of Mathematical Sciences},
  note={\url{http://dx.doi.org/10.1016/0041-5553(63)90382-3}}
}

@InProceedings{loizoustoch20,
  title = 	 { Stochastic {P}olyak Step-Size for {SGD}: An Adaptive Learning Rate for Fast Convergence },
  author =       {Loizou, Nicolas and Vaswani, Sharan and Hadj Laradji, Issam and Lacoste-Julien, Simon},
  booktitle = 	 {Proceedings of The 24th International Conference on Artificial Intelligence and Statistics},
  pages = 	 {1306--1314},
  year = 	 {2021},
  volume = 	 {130},  
  publisher =    {PMLR},  
}

@INPROCEEDINGS{sunmatrix15,
  author={Sun, Ruoyu and Luo, Zhi-Quan},
  booktitle={56th IEEE Annual Symposium on Foundations of Computer Science}, 
  title={Guaranteed Matrix Completion via Nonconvex Factorization}, 
  year={2015},
  volume={},
  number={},
  pages={270-289},
  note={\url{https://doi.org/10.1109/FOCS.2015.25}}
}

@inproceedings{Gower2021SGDFS,
  title={{SGD} for structured nonconvex functions: Learning rates, minibatching and interpolation},
  author={Gower, Robert and Sebbouh, Othmane and Loizou, Nicolas},
  booktitle={International Conference on Artificial Intelligence and Statistics},
  pages={1315--1323},
  year={2021},
  organization={PMLR}
}

@article{Bassily2018OnEC,
  title={On exponential convergence of {SGD} in non-convex over-parametrized learning},
  author={Bassily, Raef and Belkin, Mikhail and Ma, Siyuan},
  journal={arXiv preprint arXiv:1811.02564},
  year={2018},
  note={\url{https://doi.org/10.48550/arxiv.1811.02564}}
}

@article{yazdani2021asynchronous,
  title={Asynchronous parallel nonconvex optimization under the {P}olyak-{{\L}}ojasiewicz condition},
  author={Yazdani, Kasra and Hale, Matthew},
  journal={IEEE Control Systems Letters},
  volume={6},
  pages={524--529},
  year={2021},
  note={\url{https://doi.org/10.1109/LCSYS.2021.3082800}}
}

@article{hendrickson2021totally,
  author={Hendrickson, Katherine R. and Hale, Matthew T.},
  journal={IEEE Transactions on Control of Network Systems}, 
  title={Totally Asynchronous Primal-Dual Convex Optimization in Blocks}, 
  year={2023},
  volume={10},
  number={1},
  pages={454-466},
  note={\url{https://doi.org/10.1109/TCNS.2022.3203366}}
}

@INPROCEEDINGS{hendricksonhybrid21,  
author={Hendrickson, Katherine R. and Hustig-Schultz, Dawn M. and Hale, Matthew T. and Sanfelice, Ricardo G.},  
booktitle={60th IEEE Conference on Decision and Control (CDC)},   
title={Exponentially Converging Distributed Gradient Descent with Intermittent Communication via Hybrid Methods},   
year={2021},  
volume={},  
number={},  
pages={1186-1191},  
note={\url{https://doi.org/10.1109/CDC45484.2021.9683567}}
}

@book{sanfelice21,
  title={Hybrid Feedback Control},
  author={Sanfelice, Ricardo G},
  year={2021},
  publisher={Princeton University Press},
  address={Princeton, NJ, USA}
}

@book{bertsekas03,
  title={{C}onvex {A}nalysis and {O}ptimization},
  author={Bertsekas, Dimitri and Nedic, Angelia and Ozdaglar, Asuman},
  year={2003},
  publisher={Athena Scientific}
}

@INPROCEEDINGS{ubl22,
  author={Ubl, Matthew and Hale, Matthew T.},
  booktitle={61st IEEE Conference on Decision and Control (CDC)}, 
  title={Faster Asynchronous Nonconvex Block Coordinate Descent with Locally Chosen Stepsizes}, 
  year={2022},
  volume={},
  number={},
  pages={4559-4564},
  note={\url{https://doi.org/10.1109/CDC51059.2022.9993341}}
}

@inproceedings{ubl23,
  title={Linear regularizers enforce the strict saddle property},
  author={Ubl, Matthew and Hale, Matthew and Yazdani, Kasra},
  booktitle={Proceedings of the AAAI Conference on Artificial Intelligence},
  volume={37},
  pages={10017--10024},
  year={2023},
  note={\url{https://doi.org/10.1609/aaai.v37i8.26194}}
}

@article{behrendt23,
  title={A Totally Asynchronous Algorithm for Time-Varying Convex Optimization Problems},
  author={Behrendt, Gabriel and Hale, Matthew},
  journal={IFAC-PapersOnLine},
  volume={56},
  number={2},
  pages={5203--5208},
  year={2023},
  note = {\url{https://doi.org/10.1016/j.ifacol.2023.10.116}}
}

@article{wensing20,
  title={Beyond convexity—contraction and global convergence of gradient descent},
  author={Wensing, Patrick M and Slotine, Jean-Jacques},
  journal={PLOS One},
  volume={15},
  number={8},
  pages={1-29},
  year={2020},
  publisher={Public Library of Science San Francisco, CA USA},
  note={\url{https://doi.org/10.1371/journal.pone.0243330}}
}

@article{ferrante16,
  title={State estimation of linear systems in the presence of sporadic measurements},
  author={Ferrante, Francesco and Gouaisbaut, Fr{\'e}d{\'e}ric and Sanfelice, Ricardo G and Tarbouriech, Sophie},
  journal={Automatica},
  volume={73},
  pages={101--109},
  year={2016},
}
\endgroup

\begin{acknowledgements}
Research by M. T. Hale and K. R. Hendrickson partially supported by AFOSR Grants no. FA9550-19-1-0169 and FA9550-23-1-0120,
ONR Grants no. N00014-24-1-2331 and N00014-21-1-2495, and AFRL Grant no. FA8651-22-F-1052. 
Research by D. M. Hustig-Schultz and R. G. Sanfelice partially supported by NSF Grants no. CNS-2039054 and CNS-2111688, by AFOSR Grants nos. FA9550-23-1-0145, FA9550-23-1-0313, and FA9550-23-1-0678, by AFRL Grant nos. FA8651-22-1-0017 and FA8651-23-1-0004, by ARO Grant no. W911NF-20-1-0253, and by DoD Grant no. W911NF-23-1-0158.
\end{acknowledgements}

\appendix
\section{Appendix}
\subsection{Section~\ref{sec:hsm} Proofs} \label{app:hsm}
\noindent \textbf{Proof of Lemma~\ref{lem:solutions}:}\\
All gradients in~\eqref{eq:hfirst} are well-defined under Assumption~\ref{as:L}. 
Using Proposition 6.10 in~\cite{goebel12} with~$U=C$, let~$\xi=(x, \eta, \tau) \in C \backslash D$. 
Because~$C = \R^n \times \R^{nN} \times \T$, we see that the tangent cone
to~$\xi \in C \backslash D$ is 
\begin{equation}
T_C(\xi) = 
\begin{cases}
(-\infty, \infty)^n \times (-\infty, \infty)^{nN} \times (-\infty, \infty) & \tau \in (0, \tau_{\max}) \\
(-\infty, \infty)^n \times (-\infty, \infty)^{nN} \times (-\infty, 0] & \tau = \tau_{\max}
\end{cases}.
\end{equation}
Then~$f(\xi) \subset T_C(\xi)$. Because~$G(D) \subset C$, case (c) in Proposition~6.10 does not apply. We avoid case (b) of Proposition~6.10 by showing that 
there is no finite escape time for any solution. 
To that end, 
consider a maximal solution~$\phi$. Then~$\phi_{x}(0,0)$ and~$\phi_{\eta^i}(0,0)$ denote the initial values of~$\phi_{x}$ and~$\phi_{\eta^i}$ for all~$i \in [N]$, respectively. Denote the time at which agents perform their first jump as~$(t_{1},0)$. 
First, 
consider~$\phi(0, 0) \in C$. 
Then $\phi(t_1, 0) = \phi(0, 0) + t_1 f\big(\phi(0, 0)\big)$,
where~$f$ is from~\eqref{eq:hfirst}. At the first jump,~$\phi_{x_i}$
remains the same for all~$i \in [N]$ (i.e., we have
$\phi_x(t_1, 1) = \phi_x(t_1, 0)$), we set $\phi_{\eta^i} = 
\textnormal{col}(\phi_{x_1}, \ldots, \phi_{x_N})$
for all~$i \in [N]$, and~$\phi_{\tau}$ is reset
to~$[\tau_{\min}, \tau_{\max}]$,
all of which imply that~$|\phi(t_1, 1)| < \infty$. 
A similar argument proves the boundedness
of~$\phi(t_1, 1)$ when~$\phi(0, 0) \in D$.
Iterating this argument forward in time, we see that
the flow map is piecewise constant over flow intervals
and the jump map simply copies certain entries
of~$\phi_{x}$ into~$\phi_{\eta}$ in the appropriate way,
while~$\phi_{\tau}$ is always reset to a compact interval. 
Thus, repeating the preceding argument 
proves the finiteness of solutions across flow intervals
and at jump times, which rules out finite escape time. 
Therefore, $\textnormal{dom } \phi$ is unbounded and case (b)
of Proposition~6.10 in~\cite{goebel12} does not apply,
which means that case (a) of that result must hold, namely,
any maximal solution~$\phi$ is complete. 
Finally, Zeno behavior is ruled out by noting that~$\tau_{\min} > 0$.

\subsection{Section~\ref{sec:stability} Proofs} \label{app:stability}
Towards proving Lemma~\ref{lem:comparison}, we first state and prove several intermediate results. We note that given Assumption~\ref{as:L},
the function~$\nabla_i L$ is also Lipschitz, where~$\nabla_i L := \frac{\partial}{\partial x_i}L$ 
is the derivative of~$L$ with respect to the~$i^{th}$ block of~$x$.
\begin{lemma} \label{lem:lipschitz}
Let~$\nabla L$ be~$K$-Lipschitz. Then~$\nabla_i L$ is also~$K$-Lipschitz for all~$i \in [N]$, 
i.e., the inequality~$\norm{ \nabla_i L(x) - \nabla_i L(y) } \leq K \norm{x - y}$ 
holds for all~$x,y \in \R^n$.
\end{lemma}
\noindent \emph{Proof: }
From the definition of~$K$-Lipschitz, we may write~$\norm{ \nabla L(x) - \nabla L(y) } \leq K \norm{x-y}$ for all~$x,y \in \R^n$. Noting 
that~$\norm{\nabla_i L(x) - \nabla_i L(y)} \leq 
\norm{\nabla L(x) - \nabla L(y)}$ gives the desired result. \hfill $\qed$

Furthermore, based on Proposition A.32 in~\cite{bertsekas89}, the smoothness of~$L$ allows us to apply the Descent Lemma given in~Lemma~\ref{lem:descent}.

\begin{lemma}[Descent Lemma, Proposition~A.32 in~\cite{bertsekas89}] \label{lem:descent}
Let~$L : \R^n \to \R$ be continuously differentiable and have the Lipschitz property~$\norm{\nabla L(x) - \nabla L(y)} \leq K \norm{x-y}$ for every~$x,y \in \R^n$. Then 
for all~$x,y$ in~$\R^n$, 
\begin{equation}
    L(y) \leq L(x) + \nabla L(x)^{\top} (y-x) + \frac{K}{2}\norm{y-x}^2. \label{eq:descent}
\end{equation}
\end{lemma}

\noindent \textbf{Proof of Lemma~\ref{lem:comparison}:}\\
By construction,~$V(\xi)$ is zero only for~$\xi \in \mathcal{A}$ and is positive otherwise. By Assumption~\ref{as:pl}, for 
any~$x^* \in \X^*$,
we have~$L(x^*) = L^*$. For a fixed~$\xi = (x, \eta, \tau)$, we define 
\begin{equation}
    \hat{x}^0 := \argmin_{x^* \in \X^*}  \norm{x - x^*} \qquad \textnormal{ and } \qquad
    \hat{x}^i := \argmin_{x^* \in \X^*}  \norm{\eta^i - x^*} \textnormal{ for all } i \in [N].
\end{equation}
Then~$V(\xi)$ is equivalent to
    $V(\xi) = \bigl(L(x) - L(\hat{x}^0)\bigr) +  \sum_{i \in [N]} \bigl(L(\eta^i) - L(\hat{x}^i)\bigr)$.
Because~$\nabla L$ is~$K$-Lipschitz, 
Lemma~\ref{lem:descent} implies that~$L(x)-L(\hat{x}^0) \leq 
\frac{K}{2} \norm{x-\hat{x}^0}^2$ for all~$x\in \R^n$. 
For the same reason, we
have~$L(\eta^i) - L(\hat{x}^i) \leq \frac{K}{2}\norm{\eta^i - \hat{x}^i}$
for all~$\eta^i \in \R^n$.
Thus,~$V(\xi)$ may be bounded as
\begin{align}
    V(\xi) \leq \frac{K}{2} \Bigl(   \norm{x - \hat{x}^0}^2 +  \sum_{i \in [N]} \norm{\eta^i - \hat{x}^i}^2 \Bigr) = \frac{K}{2} \bigl| \xi\bigr|^2_\mathcal{A}\label{eq:lxi}.
\end{align}
Therefore, we set~$\alpha_2(s) = \frac{K}{2} s^2 \in \mathcal{K}_\infty$ for all~$s \geq 0$.

Because~$L$ is~$\beta$-PL and has a Lipschitz gradient, it also satisfies the quadratic growth condition (QG) with constant~$\beta$ (see Theorem 2 in~\cite{karimi2016linear}). In particular, given any~$x\in \R^n$, we have~$L(x) - L^* \geq \frac{\beta}{2} \min_{x^* \in \X^*} \norm{x - x^*}^2$. Thus, using the definitions of~$\hat{x}^0$ and~$\hat{x}^i$ above, we can write
$L(x) - L^* \geq \frac{\beta}{2} \norm{x - \hat{x}^0}^2$
and 
$L(\eta^i) - L^* \geq \frac{\beta}{2}  \norm{\eta^i - \hat{x}^i}^2$ 
for all~$i \in [N]$. 
This leads to
\begin{align}
    V(\xi) \geq \frac{\beta}{2} \Bigl(\norm{x - \hat{x}^0}^2 +  \sum_{i \in [N]} \norm{\eta^i - \hat{x}^i}^2 \Bigr) = \frac{\beta}{2} \bigl| \xi\bigr|^2_\mathcal{A}.
\end{align}
Setting~$\alpha_1(s) = \frac{\beta}{2} s^2 \in \mathcal{K}_\infty$ for all~$s \geq 0$ completes the proof.
\hfill $\qed$

\noindent \textbf{Proof of Proposition~\ref{prop:conv}:}\\
We first consider~$\xi \in C$ and the Lyapunov function~$V$ defined in Lemma~\ref{lem:comparison}, where
\begin{align}
    \nabla V(\xi) & = \begin{bmatrix}
    \nabla L(x) \\
    \col \bigl(\nabla L(\eta^1), \dots, \nabla L(\eta^N)\bigr) \\
    0
    \end{bmatrix}.
\end{align}
This leads to 
\begin{align}
    \bigl\langle \nabla V(\xi), f(\xi) \bigr\rangle 
    = - \nabla L(x)^{\top} h(\eta) 
    = - \sum_{i \in [N]} \nabla_i L(x)^{\top} \nabla_i L(\eta^i)
    \quad \textnormal{for all } \quad \xi \in C,\label{eq:innerdef}
\end{align}
where~$h$ is from~\eqref{eq:hdef}. 
By Lemma~\ref{lem:solutions}, we know that every maximal
solution is complete, and
we now pick a maximal solution~$\phi$ initialized such that~$\phi_{x}(0,0) = \phi_{\eta^i}(0,0)$ for all~$i \in [N]$.  For each~$I^j := \{t:(t,j) \in \textnormal{dom } \phi \}$ with nonempty interior and with~$t_{j+1} > t_j$ such that~$[t_j, t_{j+1}] = I^j$, the
initialization of~$\phi$ 
leads to a common value of~$\eta^i$ across all agents for all time, i.e.,~$\phi_{\eta^i}(t,j) = \phi_{\eta^\ell}(t,j)$ for all 
pairs of 
agents~$i$ and~$\ell$ and all times~$(t,j)$. For simplicity, we denote this shared value with~$\Bar{\eta}$. This allows us to rewrite~\eqref{eq:innerdef} as
\begin{multline}
    \bigl\langle \nabla V \bigl(\phi(t,j)\bigr), f \bigl(\phi(t,j)\bigr) \bigr\rangle  = - \sum_{i \in [N]} \nabla_i L\bigl(\phi_{x}(t,j)\bigr)^{\top} \nabla_i L\bigl(\phi_{\Bar{\eta}}(t,j)\bigr) 
    \\ = - \nabla L\bigl(\phi_{x}(t,j)\bigr)^{\top} \nabla L\bigl(\phi_{\Bar{\eta}}(t,j)\bigr). \label{eq:vdotdefined}
\end{multline}
We now apply~Lemma~\ref{lem:descent} with~$x=\phi_{x}(t,j)$ and~$y=\phi_{\Bar{\eta}}(t,j)$, giving
\begin{align}
    &L\bigl(\phi_{\Bar{\eta}}(t,j)\bigr) \leq L\bigl(\phi_{x}(t,j)\bigr) + \nabla L\bigl(\phi_{x}(t,j)\bigr)^{\top}\bigl(\phi_{\Bar{\eta}}(t,j) - \phi_{x}(t,j)\bigr) + \frac{K}{2}\norm{\phi_{\Bar{\eta}}(t,j) - \phi_{x}(t,j)}^2 \\
    &= L\bigl(\phi_{x}(t,j)\bigr) + \sum_{i \in [N]} \nabla_i L\bigl(\phi_{x}(t,j)\bigr)^{\top} \bigl(\phi_{\Bar{\eta}_i}(t,j) - \phi_{x_i}(t,j)\bigr) + \frac{K}{2}\sum_{i \in [N]} \norm{\phi_{\Bar{\eta}_i}(t,j) - \phi_{x_i}(t,j)}^2.
\end{align}
Applying the relationships~\eqref{eq:z2update} and~\eqref{eq:z1update} from Lemma~\ref{lem:contract},
we find 
\begin{align}
    L\bigl(\phi_{\Bar{\eta}}&(t,j)\bigr) \leq L\bigl(\phi_{x}(t,j)\bigr) + \frac{K}{2}\sum_{i \in [N]} \norm{\phi_{\Bar{\eta}_i}(t_j,j) - \phi_{\Bar{\eta}_i}(t_j,j) + (t-t_j) \nabla_i L\bigl(\phi_{\Bar{\eta}}(t_j,j)}^2 \\
    &\qquad+ \sum_{i \in [N]} \nabla_i L\bigl(\phi_{x}(t,j)\bigr)^{\top} \Bigl(\phi_{\Bar{\eta}_i}(t_j,j) - \phi_{\Bar{\eta}_i}(t_j,j) + (t-t_j) \nabla_i L\bigl(\phi_{\Bar{\eta}}(t_j,j)\bigr)\Bigr) \\
    &= L\bigl(\phi_{x}(t,j)\bigr) + (t-t_j) \nabla L\bigl(\phi_{x}(t,j)\bigr)^{\top} \nabla L\bigl(\phi_{\Bar{\eta}}(t,j)\bigr) + \frac{K}{2} (t - t_j)^2  \norm{\nabla L\bigl(\phi_{\Bar{\eta}}(t,j)\bigr)}^2 .
\end{align}
Rearranging gives 
\begin{multline}
    -(t-t_j) \nabla L\bigl(\phi_{x}(t,j)\bigr)^{\top} \nabla L\bigl(\phi_{\Bar{\eta}}(t,j)\bigr) \leq \\ \Bigl( L\bigl(\phi_{x}(t,j)\bigr) - L\bigl(\phi_{\Bar{\eta}}(t,j)\bigr) \Bigr) + \frac{K}{2} (t - t_j)^2  \norm{\nabla L\bigl(\phi_{\Bar{\eta}}(t,j)\bigr)}^2. \label{eq:rearranged}
\end{multline}
We now apply~\eqref{eq:descent} from Lemma~\ref{lem:descent} once more, using~$x=\phi_{\Bar{\eta}}(t,j)$ and~$y=\phi_{x}(t,j)$ to write
\begin{align}
    L\bigl(\phi_{x}(t,j)\bigr) &- L\bigl(\phi_{\Bar{\eta}}(t,j)\bigr) \leq \nabla L\bigl(\phi_{\Bar{\eta}}(t,j)\bigr)^{\top} \bigl(\phi_{x}(t,j)-\phi_{\Bar{\eta}}(t,j)\bigr) + \frac{K}{2}\norm{\phi_{x}(t,j)-\phi_{\Bar{\eta}}(t,j)}^2 \\
    &= \sum_{i \in [N]} \nabla_i L\bigl(\phi_{\Bar{\eta}}(t,j)\bigr)^{\top} \bigl(\phi_{x_i}(t,j)-\phi_{\Bar{\eta}_i}(t,j)\bigr) + \frac{K}{2} \sum_{i \in [N]} \norm{\phi_{x_i}(t,j)-\phi_{\Bar{\eta}_i}(t,j)}^2\\
    &= -(t-t_j)  \sum_{i \in [N]} \nabla_i L\bigl(\phi_{\Bar{\eta}}(t,j)\bigr)^{\top} \nabla_i L\bigl(\phi_{\Bar{\eta}}(t,j)\bigr) + \frac{K}{2} (t - t_j)^2  \norm{\nabla L\bigl(\phi_{\Bar{\eta}}(t,j)\bigr)}^2 \\
    &= -(t-t_j) \norm{\nabla L\bigl(\phi_{\Bar{\eta}}(t,j)\bigr)}^2 + \frac{K}{2} (t - t_j)^2  \norm{\nabla L\bigl(\phi_{\Bar{\eta}}(t,j)\bigr)}^2, \label{eq:smooth2}
\end{align}
where the second equality applies Lemma~\ref{lem:contract}. Applying the last inequality to~\eqref{eq:rearranged} and grouping terms, we find 
\begin{equation}
    -(t-t_j) \nabla L\bigl(\phi_{x}(t,j)\bigr)^{\top} \nabla L\bigl(\phi_{\Bar{\eta}}(t,j)\bigr) 
    \leq -(t-t_j) ( 1 - K \tau_{\max} ) 
    \norm{\nabla L\bigl(\phi_{\Bar{\eta}}(t,j)\bigr)}^2.
\end{equation}
Dividing by~$t-t_j$, which is positive by definition, gives
\begin{align}
    - \nabla L\bigl(\phi_{x}(t,j)\bigr)^{\top} \nabla L\bigl(\phi_{\Bar{\eta}}(t,j)\bigr) \leq - ( 1 - K \tau_{\max} ) 
    \norm{\nabla L\bigl(\phi_{\Bar{\eta}}(t,j)\bigr)}^2, \label{eq:innerprod0}
\end{align}
where the right hand side is negative for~$1-K\tau_{\max} > 0$, which is satisfied for~$\tau_{\max} < \frac{1}{K}$. Furthermore, the right hand side of~\eqref{eq:innerprod0} will be zero only when an optimum of~$L$ has been reached, namely at a point at which~$\nabla L$ is zero. Finally, applying the~$\beta$-PL condition from Assumption~\ref{as:pl} allows us to write
\begin{align}
    - \nabla L\bigl(\phi_{x}(t,j)\bigr)^{\top} \nabla L\bigl(\phi_{\Bar{\eta}}(t,j)\bigr) \leq - 2 \beta ( 1 - K \tau_{\max} ) \Bigl( L\bigl(\phi_{\Bar{\eta}}(t,j)\bigr) - L^* \Bigr). \label{eq:PLappl}
\end{align}
Looking once more at~\eqref{eq:smooth2}, we note that
\begin{align}
    L\bigl(\phi_{x}(t,j)\bigr) - L\bigl(\phi_{\Bar{\eta}}(t,j)\bigr) &\leq -(t-t_j) \Bigl(1 - \frac{K}{2} (t-t_j) \Bigr) 
    \norm{\nabla L\bigl(\phi_{\Bar{\eta}}(t,j)\bigr)}^2 \\
    &\leq - \tau_{\min} \Bigl(1 - \frac{K}{2} \tau_{\max} \Bigr) 
    \norm{\nabla L\bigl(\phi_{\Bar{\eta}}(t,j)\bigr)}^2, \label{eq:negjumps}
\end{align}
where the right hand side is negative since~$\tau_{\max} < \frac{1}{K}$. Thus,~$L\bigl(\phi_{x}(t,j)\bigr) \leq  L\bigl(\phi_{\Bar{\eta}}(t,j)\bigr)$ and
\begin{equation}
V \bigl(\phi(t,j)\bigr) \!=\!  \Bigl( L \bigl(\phi_{x}(t,j)\bigr) - L^* \Bigr) \!+\!\!  \sum_{i \in [N]} \Bigl(L\bigl(\phi_{\eta^i}(t,j)\bigr) - L^* \Bigr) 
\leq (N+1) \Bigl( L\bigl(\phi_{\Bar{\eta}}(t,j)\bigr) - L^*  \Bigr).
\end{equation}
Combining this inequality with~\eqref{eq:PLappl}, we have
\begin{align}
    - \nabla L\bigl(\phi_{x}(t,j)\bigr)^{\top} \nabla L\bigl(\phi_{\Bar{\eta}}(t,j)\bigr) &\leq - \frac{2}{N+1} \beta ( 1 - K \tau_{\max} ) V \bigl(\phi(t,j)\bigr) \label{eq:tinbound}.
\end{align}
Combined with~\eqref{eq:vdotdefined}, we conclude that
\begin{equation} \label{eq:dotprodbound}
\left\langle V(\phi(t,j)), f(\phi(t,j)) \right\rangle \leq - \frac{2}{N+1} \beta ( 1 - K \tau_{\max} ) V \bigl(\phi(t,j)\bigr) \leq 0. 
\end{equation}
%
%
%
%

Next consider the change of~$V$ at jumps. For a maximal solution~$\phi$ such that~$\phi_{x}(0,0) = \phi_{\eta^i}(0,0)$ for all~$i \in [N]$, we may write the change at jump~$j+1$ as
\begin{align}
    V\bigl(&G(\phi(t_{j+1},j)\bigr)\bigr) \!-\! V\bigl(\phi(t_{j+1},j)\bigr) = \Bigl(\! L\bigl(\phi_x(t_{j+1},j \!+\! 1)\bigr) \!-\! L^* \!\Bigr) \!+\! N\Bigl(\! L\bigl(\phi_{\Bar{\eta}}(t_{j+1},j \!+\! 1)\bigr) - L^* \!\Bigr)  \\
    &\qquad\qquad\qquad\qquad\qquad -\Bigl(\! L\bigl(\phi_x(t_{j+1},j)\bigr) - L^* \!\Bigr) -  N \Bigl(\! L\bigl(\phi_{\Bar{\eta}}(t_{j+1},j)\bigr) - L^* \!\Bigr) \\
    &=  N \Bigl( L\bigl(\phi_{\Bar{\eta}}(t_{j+1},j+1)\bigr) -   L\bigl(\phi_{\Bar{\eta}}(t_{j+1},j)\bigr) \Bigr) 
    = N \Bigl( L\bigl(\phi_{x}(t_{j+1},j)\bigr) - L\bigl(\phi_{\Bar{\eta}}(t_{j+1},j)\bigr) \Bigr)
\end{align}
for all~$(t_{j+1},j), (t_{j+1},j+1) \in \textnormal{dom }\phi$. For this quantity to be nonpositive, it is sufficient to show that~$L\bigl(\phi_{x}(t_{j+1},j)\bigr) \leq L\bigl(\phi_{\Bar{\eta}}(t_{j+1},j)\bigr)$. This inequality 
follows directly from~\eqref{eq:negjumps} by setting~$t = t_{j+1}$, where
\begin{align}
L\bigl(\phi_{x}(t_{j+1},j)\bigr) - L\bigl(\phi_{\Bar{\eta}}(t_{j+1},j)\bigr) &\leq - \tau_{\min} \Bigl(1 - \frac{K}{2} \tau_{\max} \Bigr) 
\norm{\nabla L\bigl(\phi_{\Bar{\eta}}(t_{j+1},j)\bigr)}^2.
\end{align} 
For~$\tau_{\max} < \frac{1}{K}$, the term~$1 - \frac{K}{2} \tau_{\max}$ is positive. Then~$L\bigl(\phi_{x}(t_{j+1},j)\bigr) \leq L\bigl(\phi_{\Bar{\eta}}(t_{j+1},j)\bigr)$ and 
\begin{align}
    V\bigl(G(\phi(t_{j+1},j)\bigr)\bigr) &- V \bigl(\phi(t_{j+1},j)\bigr) \leq 0. \label{eq:jinbound}
\end{align}

Following the work done in~\cite{chai17} and~\cite{goebel12}, we are able to perform direct integration in order to upper bound~$V \bigl(\phi(t,j)\bigr)$ in terms of~$V \bigl(\phi(0,0)\bigr)$ using~\eqref{eq:dotprodbound} and~\eqref{eq:jinbound} as bounds. Thus,
\begin{align}
    V \bigl(\phi(t,j)\bigr) \leq \exp \Big( - \frac{2}{N+1} \beta ( 1 - K \tau_{\max}) t\Big) V \bigl(\phi(0,0)\bigr).
\end{align}
Using the comparison functions given in Lemma~\ref{lem:comparison}, we find a bound for~$\bigl|\phi(t,j)\bigr|^2$ via
\begin{align}
    \bigl|\phi(t,j)\bigr|^2_\mathcal{A} &\leq \frac{2}{\beta} \exp \Big( - \frac{2}{N+1} \beta ( 1 - K \tau_{\max}) t\Big) V \bigl(\phi(0,0)\bigr) \\
    &\leq \frac{K}{\beta} \exp \Big( - \frac{2}{N+1} \beta ( 1 - K \tau_{\max}) t\Big) \bigl|\phi(0,0)\bigr|^2_\mathcal{A},
\end{align}
where taking the square root gives the final bound. \hfill $\qed$

\noindent \textbf{Proof of Proposition~\ref{prop:ges}: }\\
Two initialization scenarios must be considered: $\phi_{x}(0,0) = \phi_{\eta^i}(0,0)$ for all~$i \in [N]$ and $\phi_{x}(0,0) \neq \phi_{\eta^i}(0,0)$ for at least one~$i \in [N]$. For the first case,~$\phi_{x}(0,0) = \phi_{\eta^i}(0,0)$ for all~$i \in [N]$, Proposition~\ref{prop:conv} applies in its original form. This is the best-case scenario that results in the smallest upper bound on the distance to~$\mathcal{A}$. 

Now consider the second case, when~$\phi_{x}(0,0) \neq \phi_{\eta^i}(0,0)$ for at least one~$i \in [N]$. After agents have completed at least one jump, all assumptions of Proposition~$\ref{prop:conv}$ hold. Denote the time at which the first jump occurs as~$(t_1, 1)$. Thus, for any maximal solution~$\phi$ 
(which is complete by Lemma~\ref{lem:solutions}), for any
$(t,j) \in \textnormal{dom } \phi$ such that~$j\geq 1$, the following holds: 
\begin{align}
    \bigl|\phi(t,j)\bigr|_\mathcal{A} &\leq \sqrt{\frac{K}{\beta}} \exp \Big( - \frac{\beta}{N+1} ( 1 - K \tau_{\max}) t\Big) \bigl|\phi(t_{1},1)\bigr|_\mathcal{A}. \label{eq:jgeq2}
\end{align}

We now seek to bound~$\bigl|\phi(t_{1},1)\bigr|_\mathcal{A}$ in terms of~$\bigl|\phi(0,0)\bigr|_\mathcal{A}$. 
We first 
define the point $\bar{x}^0 := \argmin_{x^* \in \X^*}  \norm{\phi_{x} (t_{1}, 1) - x^*}$, and, for each~$i \in [N]$, we define~$\bar{x}^i := \argmin_{x^* \in \X^*} \norm{\phi_{\eta^i}(t_1, 1) - x^*}$.  
We begin by expanding~$\bigl|\phi(t_{1},1)\bigr|^2_\mathcal{A}$ to find 
\begin{align}
    \bigl|\phi(t_{1},1)\bigr|^2_\mathcal{A} &= 
    \norm{\phi_{x} (t_{1}, 1) - \bar{x}^0}^2 
    + \sum_{i \in [N]} \norm{\phi_{\eta^i}(t_{1}, 1) - \bar{x}^i}^2. \label{eq:distance_exp2}
\end{align}
We now define~$\hat{x}^0 := \argmin_{x^* \in \X^*}  
\norm{\phi_{x}(0,0) - x^*}$. Note that by definition of~$\bar{x}^0$ and~$\bar{x}^i$, the following inequalities hold:
\begin{align}
   \norm{\phi_{x} (t_{1}, 1) - \bar{x}^0}^2 &\leq \norm{\phi_{x} (t_{1}, 1) - \hat{x}^0}^2 = \sum_{i \in [N]} \norm{\phi_{x_i} (t_{1}, 1) - \hat{x}^0_i}^2  \\
    \norm{\phi_{\eta^i}(t_{1}, 1) - \bar{x}^i}^2 &\leq \norm{\phi_{\eta^i} (t_{1}, 1) - \hat{x}^0}^2.
\end{align}
These inequalities allow us to rewrite~\eqref{eq:distance_exp2} as
\begin{align}
\bigl|\phi(t_{1},1)\bigr|^2_\mathcal{A} &\leq  
\sum_{i \in [N]} \norm{\phi_{x_i} (t_{1}, 1) - \hat{x}^0_i}^2 + \sum_{i \in [N]} \norm{\phi_{\eta^i}(t_{1}, 1) - \hat{x}^0}^2.\label{eq:distrewrite2}
\end{align}
Define~$\hat{x}^i := \argmin_{x^* \in \X^*} 
\norm{\phi_{\eta^i}(0,0) - x^*}$ for all~$i \in [N]$. We first upper bound~$\norm{\phi_{x_i} (t_{1}, 1) - \hat{x}^0_i}^2$ by applying Lemma~\ref{lem:contract} and using~$\norm{a-b}^2 
\leq 2\norm{a}^2 + 2\norm{b}^2$, resulting in
\begin{align}
    \norm{\phi_{x_i} (t_{1}, 1) - \hat{x}^{0}_i}^2 &= 
    \norm{\phi_{x_i} (t_{1}, 0) - \hat{x}^{0}_i}^2
    = \norm{\phi_{x_i} (0,0) - t_{1} \nabla_i L\bigl(\phi_{\eta^i} (0,0)\bigr) - \hat{x}^{0}_i}^2 \\
    &\leq 2 \norm{\phi_{x_i} (0,0) -  \hat{x}^{0}_i}^2 + 2 K^2 \tau_{\max}^2 \norm{\phi_{\eta^i}(0,0) - \hat{x}^i}^2, \label{eq:xupperPL1}
\end{align}
where the first equality applies~\eqref{eq:hlast}, the second equality applies Lemma~\ref{lem:contract}, the first inequality uses~$\nabla L(\hat{x}^i) = 0$, and the final inequality applies Lemma~\ref{lem:lipschitz}.
Summing over all~$i$ on both sides of~\eqref{eq:xupperPL1} gives
\begin{equation}
    \norm{\phi_{x} (t_{1}, 1) - \hat{x}^{0}}^2 
    \leq 2 \norm{\phi_{x} (0,0) -  \hat{x}^{0}}^2 + 2 K^2 \tau_{\max}^2 \sum_{i \in [N]} \norm{\phi_{\eta^i} (0,0) -  \hat{x}^i}^2. \label{eq:xupperPL2}
\end{equation}

An upper bound on~$\norm{\phi_{\eta^i}(t_{1}, 1) - \hat{x}^0}^2$ also needs to be derived to upper bound the right hand side of~\eqref{eq:distrewrite2}. Recall that at any jump~$j$, we have~$\phi_{\eta^i}(t_j, j) = \phi_{x}(t_j, j)$ for all~$i$. Thus,~$\phi_{\eta^i}(t_1, 1) = \phi_{x}(t_1, 1)$ for all~$i$. We use this to expand and upper bound the following:
\begin{align}
    \sum_{i \in [N]} \norm{\phi_{\eta^i}(t_{1}, 1) - \hat{x}^0}^2 &= 
    N \norm{\phi_{x}(t_{1}, 1) - \hat{x}^0}^2 \\
    &\leq 2N \norm{\phi_{x} (0,0) -  \hat{x}^0}^2 + 2 N K^2 \tau_{\max}^2 \sum_{i \in [N]} \norm{\phi_{\eta^i} (0,0) -  \hat{x}^i}^2, \label{eq:etaupperPL1}
\end{align}
where the inequality applies~\eqref{eq:xupperPL2}. Summing~\eqref{eq:xupperPL2} and~\eqref{eq:etaupperPL1} gives
\begin{align}
    \norm{\phi_{x}&(t_{1}, 1) - \hat{x}^0}^2 
    + \sum_{i \in [N]} \norm{\phi_{\eta^i}(t_{1}, 1) - \hat{x}^0}^2 \leq 2 \norm{\phi_{x} (0,0) -  \hat{x}^0}^2 \\
    &+ 2 K^2 \tau_{\max}^2 \sum_{i \in [N]} \norm{\phi_{\eta^i} (0,0) -  \hat{x}^i}^2 
    \!\! + 2N \norm{\phi_{x} (0,0) -  \hat{x}^0}^2 \!\! + \! 2 N K^2 \tau_{\max}^2 \sum_{i \in [N]} \norm{\phi_{\eta^i} (0,0) -  \hat{x}^i}^2 \\
    &= 2(N\!+\!1) \norm{\phi_{x} (0,0) -  \hat{x}^0}^2 + 2 (N\!+\!1) K^2 \tau_{\max}^2 \sum_{i \in [N]} \norm{\phi_{\eta^i} (0,0) -  \hat{x}^i}^2 \\
    &\leq 2(N\!+\!1) \Bigl( \norm{\phi_{x} (0,0) -  \hat{x}^0}^2 + \sum_{i \in [N]} \norm{\phi_{\eta^i} (0,0) -\hat{x}^i}^2\Bigr),
\end{align}
where the first equality groups terms and the second inequality uses~$\tau_{\max} < \frac{1}{K}$ to simplify. Applying~\eqref{eq:distrewrite2} on the left-hand side and the definition of~$|\cdot|^2_\mathcal{A}$ on the right-hand side gives
\begin{align}
    \bigl|\phi(t_1,1)\bigr|^2_\mathcal{A} &\leq 2(N + 1) \bigl|\phi(0,0)\bigr|^2_\mathcal{A}. \label{eq:phit1bound}
\end{align}
Taking the square root and applying the resulting inequality to~\eqref{eq:jgeq2} gives a bound for all~$j \geq 1$. 
\hfill $\qed$

\end{document}